\documentclass[a4paper,12pt]{amsart}

\usepackage[
  margin=30mm,
  marginparwidth=25mm,     
  marginparsep=2mm,       
  bottom=25mm,
  ]{geometry}

\usepackage[T1]{fontenc}
\usepackage[utf8]{inputenc}
\usepackage[english]{babel}
\usepackage{amsmath,amssymb,amsthm,mathtools}
\usepackage{latexsym}
\usepackage{delarray}
\usepackage{bbm}
\usepackage{scrextend}
\usepackage{hyperref}
\usepackage[pdftex,usenames,dvipsnames]{xcolor}
\usepackage{datetime}
\usepackage{mathrsfs}
\usepackage{enumitem}  
\usepackage{tikz}
\usepackage{comment}
\usepackage{multicol}
\usepackage[makeroom]{cancel}

\newtheorem{Th}{Theorem}[section]

\newtheorem{Lem}[Th]{Lemma}

\newtheorem{Rem}[Th]{Remark}

\newcommand{\R}{\mathbb{R}}

\newcommand{\C}{\mathbb{C}}

\usepackage{isomath}

\DeclareMathOperator{\sgn}{sgn}

\allowdisplaybreaks 

\title[A note on the quartic gKdV equation in weighted Sobolev spaces]
{A note on the quartic generalized Korteweg-de Vries equation in weighted Sobolev spaces}

\author[A. J. Castro]{Alejandro J. Castro}
\author[A. Esfahani]{Amin Esfahani}
\author[L. Zhapsarbayeva]{Lyailya Zhapsarbayeva}

\address{\newline        
Alejandro J. Castro,  
Amin Esfahani,
Lyailya Zhapsarbayeva \newline        
Department of  Mathematics, Nazarbayev University, \newline 010000 Astana, Kazakhstan}
\email{\{alejandro.castilla,amin.esfahani,lyailya.zhapsarbayeva\}@nu.edu.kz}

 \thanks{
This paper was supported by the Ministry of Education and Science of the Republic of Kazakhstan (Grant No AP19676408)}



\begin{document}
\maketitle
\vspace*{-0.5cm}
\begin{abstract}
In this paper we establish the persistence property for solutions of the quartic 
generalized Korteweg-de Vries equation with initial data in  weighted Sobolev spaces  
$H^{s}(\R)\cap L^2(|x|^{2r}dx)$ for 
$s =1/12 + \varepsilon$
and any $r\in (0,R)$, for some $0 < \varepsilon < 1/4$ and $0<R<s/2$. 
\end{abstract}

\section{Introduction}
In this work we are interested in the initial value problem for the k-generalized Korteweg-de Vries equation (k-gKdV)
\begin{equation}\label{eq:k-gKdV}
\left\{
\begin{array}{ll}
\partial_t u +\partial_x^3 u
+  \partial_x (u^{k+1})  =0, & x \in \mathbb{R}, \, t>0, \\
u(x,0)=u_0(x), & x \in \mathbb{R},
\end{array}
\right.
\end{equation}
which was originally proposed to model the unidirectional propagation of nonlinear dispersive water waves (\cite{KdV}), but has also been considered in connection with other physical systems, such as elastic rods (\cite{Nar}), gravity waves (\cite{Li2009})  and plasma physics (\cite{BeKa}).\\

Next, we briefly review well-posedness results
(i.e. existence, uniqueness and stability of solutions) obtained for the k-gKdV equation with initial data in Sobolev spaces $H^s(\R)$, with emphasis on the minimal regularity $s \in \R$ possible.\\

When $k=1$,  \eqref{eq:k-gKdV} is simply referred as the Korteweg-de Vries equation (KdV). 
Z. Guo  \cite{Guo} and N. Kishimoto  \cite{Kis2009} showed independently that the KdV  is locally and globally (in time)  well-possed for initial data in $H^{-3/4}(\mathbb{R})$. \\

For $k=2$, \eqref{eq:k-gKdV} is known as the  modified Korteweg-de Vries equation (mKdV).
Its local well-posedness, provided that $s \geq 1/4$, was shown  by C. Kenig, G. Ponce and L. Vega \cite{KPV1993} and extended globally by J. Colliander, M. Keel, G. Staffilani, H. Takaoka and T. Tao \cite{CKSTT} (see also \cite{Guo,Kis2009} for the endpoint $s=1/4$).\\ 

If $k=3$,
the best results for $s \geq -1/6$ are due to
A. Gr\"unrock \cite{Grun2005} and T. Tao \cite{Tao2007}.\\

Finally, for $k \geq 4$ and $s >(k-4)/(2k)$ the local well-posedness was established in \cite{KPV1993}. However, even for $k=4$, blow up in finite time  has been shown by Y. Martel and F. Merle \cite{MarMer}.
\\

Building upon the seminal work of T. Kato \cite{Kato1983}, the solvability of various dispersive nonlinear equations has been explored within the framework of weighted Sobolev spaces $H^{s}(\R)\cap L^2(|x|^{2r}dx)$, aiming to achieve enhanced control over the decay at infinity of solutions (as evidenced in studies such as \cite{
BusJi2018,
CN2015,
CJZ2022,
FonPn,
NP2009}
and their associated references).
Notably, for equations like the modified Korteweg-de Vries (mKdV) and the generalized Korteweg-de Vries (k-gKdV) with $k \geq 4$, optimal results within this context were established by J. Nahas \cite{Nah2012}. Moreover, A. Mu\~noz-Garc\'ia \cite{MG2019} established the local well-posedness of the Korteweg-de Vries equation for $s>3/4$ in weighted Sobolev spaces. In the same spirit, our objective is to extend these investigations to encompass the quartic generalized Korteweg-de Vries equation (in the case of $k=3$), as elaborated upon in the subsequent sections.\\

As   is standard, we are going to apply the contraction mapping principle to the integral equation version of the IVP  \eqref{eq:k-gKdV} with $k=3$, i.e. 
\begin{equation}\label{InEKdV}
u
=
e^{-t\partial^3_{x}}u_0
-\int^t_{0}
e^{-(t-t')
\partial^3_{x}}\partial_x (u^4) (\cdot,t')\, dt', 
\end{equation}
where $e^{-t\partial^3_{x}}$ is the Airy semigroup introduced in Subsection \ref{Airy}. \\

We now state our main result.

\begin{Th}\label{Th:Main} 
Let 
$u_{0}\in H^{1/12+\varepsilon}(\R)\cap L^2(|x|^{2r}dx)$
and 
$r\in (0,R)$, for some
$0<\varepsilon <1/4$
and certain $0<R <1/24+\varepsilon/2$.
Then, there exist $T>0$ and a unique 
 solution $u$ of the integral equation \eqref{InEKdV} such that 
$$u(\cdot,t)\in H^{1/12+\varepsilon}(\R)\cap L^2(|x|^{2r}dx),
\quad t \in (0,T].$$ 
\end{Th}

In light of the Sobolev optimal results mentioned earlier for $k=3$, it is natural to expect an improvement of Theorem \ref{Th:Main} on the Sobolev regularity $s$, at least for some $0<s<1/12$. Similar observation applies in the situation $k=1$, if one compares \cite{Guo,Kis2009} with \cite{MG2019}.
We believe that a treatment of the  factor $|x|^r u$ 
(see Section \ref{Subsec:weightednorm} below) in the realm of Bourgain-type  spaces might be advantageous, as has been demonstrated in the low regularity unweighted cases 
(\cite{Bour93II,Grun2005,Guo,KPVCP,Kis2009,Tao2007} among others).
We are currently investigating this approach.

\section{Definitions and Preliminaries}

In this article we write $A \lesssim B$, if there exists a constant $C>0$ such that $A \leq C B$. Moreover, $A \sim B$ represents $A \lesssim B$ and $B \lesssim A$. 

\subsection{Function spaces}
   
For a function $f\in L^2(\R)$, consider its Fourier transform 
$$ 
 \widehat{f} (\xi)
 :=  
 \int_{\R} e^{-ix\xi} f(x)
\, dx, \quad \xi\in \R
$$
and its inverse Fourier transform  by
$$ 
f^{\vee}(x)
:= \frac{1}{2\pi} 
\int_{\R} e^{ix\xi} f(\xi)
\, d\xi, \quad x\in \R. 
$$

\quad \\
In this work we use the inhomogeneous Sobolev space $H^s(\mathbb{R})$, of order $s \in \mathbb{R}$, defined via the norm
\begin{equation*}
\|f\|_{H^s}
:=
\Big(  \int_{\R}
(1+|\xi|^2)^s
|\widehat{f}(\xi)|^2 \, d\xi\Big)^{1/2},
\end{equation*}
which satisfy the inclusion $H^{s'}(\mathbb{R}) \subset H^{s}(\mathbb{R})$ for $s \leq s'$, that is,
\begin{equation}\label{eq:Sobolevinclusion}
\|f\|_{H^s}
\lesssim
\|f\|_{H^{s'}}.
\end{equation}

\quad \\
In order to measure the regularity of functions defined in the space--time domain $\R \times [0,T]$ we introduce the  mixed--norm Lebesgue spaces $L^p_xL^q_T$ or $L^q_TL^p_x$, $1 \leq p,q \leq \infty$, given respectively by the norms
$$
\|f\|_{L^p_xL^q_T}
:=\Big\{ 
\int_{\R}
\Big(
\int^T_0 |f(x,t)|^q \, dt
\Big)^{p/q} \, dx
\Big\}^{1/p}  
$$
and
$$
\|f\|_{L^q_TL^p_x}
:=
\Big\{ 
\int^T_0
\Big(
\int_{\R}  |f(x,t)|^p \, dx
\Big)^{q/p} \, dt
\Big\}^{1/q}
$$
with the standard modifications involving the essential supremum when $p$ or $q$ are equal to infinity.

\subsection{Fractional derivatives} 

For $\alpha \in \C$, we define the fractional derivative $D_x^{\alpha}$ as the Fourier multiplier given by
$$(D_x^{\alpha}f)^{\wedge}(\xi)
:=
|\xi|^{\alpha}\widehat{f}(\xi).
$$
Analogously, we introduce the operator $(1+D_x^2)^\alpha$ via
$$
\big((1+D_x^2)^{\alpha}f\big)^{\wedge}(\xi)
:=
(1+|\xi|^2)^{\alpha}\widehat{f}(\xi).
$$
Hence, the Plancherel  identity allows us to write
\begin{equation*}
\|f\|_{H^s}
\sim 
\|(1+D_x^2)^{s/2}f\|_{L^2}
\lesssim 
\|f\|_{L^2}
+
\|D_x^s f\|_{L^2}.
\end{equation*}
Also, if we invoke the Hilbert transform $\mathcal H$ determined by
\begin{equation*}
(\mathcal H f)^{\wedge}(\xi)
:= - i \sgn(\xi) \widehat{f}(\xi),
\end{equation*}
we can relate $D_x$ with the standard derivative $\partial_x$ as $D_x=\mathcal H \partial_x$ 
or $\partial_x=\mathcal H D_x$.\\

A   highly useful property in Section \ref{Sect:ProofTh} will be the following fractional Leibniz rule-type inequality (see \cite[Theorem A.8]{KPV1993})
\begin{equation}\label{thA8}
\|
D_x^\alpha(fg)
- f D_x^\alpha g
- g D_x^\alpha f
\|_{L^p_x L^q_T}    
\lesssim 
\|D_x^{\alpha_1} f\|_{L^{p_1}_x L^{q_1}_T}
\|D_x^{\alpha_2} g\|_{L^{p_2}_x L^{q_2}_T},
\end{equation}
where $\alpha \in (0,1)$, $\alpha_1, \alpha_2 \in [0,\alpha]$ verify $\alpha = \alpha_1 + \alpha_2$
and the exponents $p, p_1, p_2, q, q_1, q_2 \in (1,\infty)$
are given by
\begin{equation}\label{eq:exponents}
\frac{1}{p} = \frac{1}{p_1} + \frac{1}{p_2}, \quad 
\frac{1}{q} = \frac{1}{q_1} + \frac{1}{q_2}.
\end{equation}
Moreover, for $\alpha_1=0$ the value $q_1=\infty$ is allowed.

\quad \\
We also recall the chain rule for fractional derivatives (\cite[Theorem A.6]{KPV1993})
\begin{equation}\label{Sta4.6}
\|D_x^{\alpha} F(f)\|_{L^{p}_x L^{q}_T}
\lesssim
\|F'(f)\|_{L^{p_1}_x L^{q_1}_T}
\|D_x^{\alpha} f\|_{L^{p_2}_x L^{q_2}_T}
\end{equation}
for $\alpha \in (0,1)$, $p, p_1, p_2, q, q_2 \in (1,\infty)$
and $q_1 \in (1,\infty]$
satisfying \eqref{eq:exponents}. \\

Sometimes it is more convenient to use this other definition of fractional derivative, 
\begin{equation*}
\mathscr{D}^{\alpha}_xf(x)
:=
\Big(
\int_{\R}\frac{|f(x)-f(y)|^2}{|x-y|^{1+2\alpha}} \, dy
\Big)^{1/2}, \quad 
\alpha \in (0,1),
\end{equation*}
which is related to $D_x^\alpha$ in the following sense (see for example \cite[eq. (1.23)]{NP2009})
\begin{equation}\label{eq:5}
\|f\|_{L^{p}}
+\|D^{\alpha}_x f\|_{L^{p}}
\sim  
\|f\|_{L^{p}}
+\|\mathscr{D}^{\alpha}_x f \|_{L^{p}},
\end{equation}
 and also satisfies (\cite[eq. (2.2)]{NP2009})
\begin{equation}\label{eq:6}
\|\mathscr{D}^{\alpha}_x(fg)\|_{L^{p}}
\lesssim 
\|f\|_{L^{\infty}} \|\mathscr{D}^{\alpha}_xg\|_{L^{p}}
+
\|g\mathscr{D}^{\alpha}_xf\|_{L^{p}}.
\end{equation}

\subsection{Airy semigroup}\label{Airy}
The solution of the Airy equation
\begin{equation*}
\left\{
\begin{array}{ll}
\partial_t u + \partial_x^3 u=0, & x \in \mathbb{R}, \, t>0, \\
u(x,0)=u_0(x), & x \in \mathbb{R}
\end{array}
\right.
\end{equation*}
can be written as $u(x,t)=e^{-t\partial^3_{x}}u_0(x)$, where $e^{-t\partial^3_{x}}$ denotes the Fourier multiplier given by
\begin{equation*}
(e^{-t\partial^3_{x}} u_0)^{\wedge}(\xi)
:= e^{it\xi^3} \widehat{u_0}(\xi).
\end{equation*}

\quad \\
The Plancherel theorem easily provides
\begin{equation}\label{2.5lem2.4}
\| e^{-t\partial^3_{x}}u_0\|_{L^{2}_x}
\sim 
\|u_0 \|_{L^{2}}.    
\end{equation}
Next, we recall additional boundedness properties of the family of operators $\{e^{-t\partial^3_{x}}\}_{t>0}$, which will be very helpful in proving our main result.
In \cite[Theorem 3.5]{KPV1993} it was established that
\begin{equation}\label{ineq36}
\|\partial_x e^{-t\partial^3_{x}}u_0 \|_{L^{\infty}_{x}L^{2}_T}
\lesssim \|u_0 \|_{L^{2}}. 
\end{equation} 

In addition, for any
$\theta \in [0,1]$, 
$\alpha \in [0,1/2]$, 
$p=2/(1-\theta)$
and 
$q=6/(\theta\alpha+\theta)$
it was shown in \cite[Lemma 2.4]{KPV1989} that
\begin{equation}\label{eq:KPV1989}
\|D^{\theta\alpha/2}_x 
e^{-t\partial^3_{x}}u_0\|_{L^q_T L_x^p}
\lesssim
\|u_0\|_{L^2_x}.
\end{equation}
and 
\begin{equation*}
\Big\| 
\int^t_0 D^{\theta\alpha}_xe^{-(t-t')\partial^3_{x}} f(\cdot,t')\, dt'
\Big\|_{L^{q}_T L^{p}_x}\leq \|f\|_{L^{q'}_TL^{p'}_x},
\end{equation*}
provided   
$1/p + 1/p'
=
1
=
1/q + 1/q'$. \\


We continue with some other crucial estimates involving integrals of the semigroup. In 
 \cite[Theorem 3.5(ii)]{KPV1993} it was proved that
\begin{equation}\label{nah3.3}
\Big\|
\partial_x \int^t_0 
e^{-(t-t')\partial^3_{x}}f(\cdot,t')\, dt'
\Big\|_{L^\infty_T L^2_x}
\lesssim
\|f\|_{L^{1}_xL^{2}_T}  
\end{equation}
and 
\begin{equation}\label{kpv-es-3}
\Big\|
\partial_x^2 \int^t_0 e^{-(t-t')\partial^3_{x}}f(\cdot,t')\, dt'
\Big\|_{ L^\infty_xL^2_T}
\lesssim
\|f\|_{L^{1}_xL^{2}_T} , 
\end{equation}
and 
in \cite[eq. (2.13)]{GnTst89}
\begin{equation}\label{eq:GinebreTsutsumi1}
\Big\| 
\int^t_0 e^{-(t-t')\partial^3_{x}} f(\cdot,t')\, dt'
\Big\|_{L^{q_1}_T L^{p_1}_x}\leq \|f\|_{L^{q'_2}_TL^{p'_2}_x},
\end{equation}
provided   
$1/p_2 + 1/p'_2
=
1
=
1/q_2 + 1/q'_2$
and
\begin{equation*}
\frac{1}{q_j}=\frac{1}{6}-\frac{1}{3p_j}, \quad 
2 \leq p_j \leq \infty
, \quad
j=1,2.
\end{equation*}


Next, we deduce other related estimates that will be used in the proof of Theorem \ref{Th:Main}. 
\begin{Lem}\label{PrL24L8/3}
For $0<\varepsilon <1/4$, one has that
\begin{equation}\label{eq:mainLem2.1}
\|  e^{-t\partial^3_{x}}u_0\|_{L_x^{120/(26+51\varepsilon)}L_T^{30/(2-3\varepsilon)}}
\lesssim \|D^{1/12+\varepsilon}_x u_0\|_{L^{2}}
\end{equation}
and
\begin{equation}\label{eq:main2Lem2.1}
\|  e^{-t\partial^3_{x}}u_0\|_{L_x^{60/(13-12\varepsilon)}L_T^{15/(1-4\varepsilon)}}
\lesssim \|u_0\|_{H^{1/12+\varepsilon}}.
\end{equation}
\end{Lem}

\begin{proof}
It follows from 
Stein's theorem of analytic interpolation
(see \cite[Theorem 1, p. 313]{BenPan}) to the family of operators 
$$
T^{(1)}_z
:=
D_x^{-\frac{2z}{3}-\frac{1}{12}}e^{-t\partial^3_{x}} 
$$
for \eqref{eq:mainLem2.1}  and 
$$
T^{(2)}_z
:=
(1+D^2_x)^{-\frac{z}{8}-\frac{1}{24}}e^{-t\partial^3_{x}} 
$$
for \eqref{eq:main2Lem2.1}, 
where $z\in \C$, with $0\leq \textrm{Re}\, z\leq 1$.
More precisely, for any $\gamma \in \R$, observe  from
\cite[Lemma 3.26]{KPV1993}  that
\begin{equation*}
\| 
T^{(1)}_{i \gamma}u_0
\|_{L^{60/13}_xL^{15}_T}
=
\|
D^{-1/12}_x 
e^{-t\partial^3_{x}}
D_x^{-2i \gamma/3}
u_0\|_{L^{60/13}_xL^{15}_T}
\lesssim 
\|D_x^{-2i \gamma/3} u_0 \|_{L^2}
\sim 
\|u_0 \|_{L^2}.
\end{equation*}
Moreover, we obtain from \cite[eq. (3.39.a)]{KPV1993} that
\begin{equation*}
 \| T^{(1)}_{1+i\gamma}u_0\|_{L^{2}_x L^{\infty}_T}
 =
 \|D^{-3/4}_x
 e^{-t\partial^3_{x}}
 D_x^{-2i \gamma/3} u_0 \|_{L^{2}_x L^{\infty}_T}
 \lesssim \|u_0 \|_{L^2}.    
\end{equation*}

Similarly, 
\cite[Lemma 3.26]{KPV1993} implies that
\begin{equation*}
\| 
T^{(2)}_{i \gamma}u_0
\|_{L^{60/13}_xL^{15}_T}
=
\|(1+D^2_x)^{-\frac{1}{24}}e^{-t\partial^3_{x}}
(1+D^2_x)^{-\frac{i\gamma}{8}}
u_0\|_{L^{60/13}_xL^{15}_T}
\lesssim 
\|u_0 \|_{L^2}
\end{equation*}
and 
 \cite[Theorem A.]{Sjo2011} gives 
\begin{equation*}
 \| T^{(2)}_{1+i\gamma}u_0\|_{L^{6}_x L^{\infty}_T}
 =
 \|(1+D^2_x)^{-\frac{1}{6}}e^{-t\partial^3_{x}}
(1+D^2_x)^{-\frac{i\gamma}{8}}
u_0 \|_{L^{6}_x L^{\infty}_T}
 \lesssim \|u_0 \|_{L^2}  .  
\end{equation*}
Therefore, we have for any $0<\theta_1,\theta_2 < 1$ and $1 \leq p_1,p_2,q_1,q_2 \leq \infty$, satisfying
\begin{equation*}
\frac{1}{p_1}
= 
\frac{1-\theta_1}{\frac{60}{13}} + \frac{\theta_1}{2}, \quad
\frac{1}{q_1}
= \frac{1-\theta_1}{15} + \frac{\theta_1}{\infty} 
\end{equation*}
and
\begin{equation*}
\frac{1}{p_2}
= 
\frac{1-\theta_2}{\frac{60}{13}} + \frac{\theta_2}{6}, \quad
\frac{1}{q_2}
= \frac{1-\theta_2}{15} + \frac{\theta_2}{\infty},
\end{equation*}
 that
\begin{equation*}
\|T^{(1)}_{\theta_1} u_0\|_{L^{p_1}_{x}L^{q_1}_T}  
\lesssim
\| u_0\|_{L^2} 
\end{equation*}
 and
\begin{equation*}
\|T^{(2)}_{\theta_2} u_0\|_{L^{p_2}_{x}L^{q_2}_T}  
\lesssim
\| u_0\|_{L^2}. 
\end{equation*}
Thus, \eqref{eq:mainLem2.1} and \eqref{eq:main2Lem2.1} are derived by choosing $\theta_1=3\varepsilon/2$ and $\theta_2=4\varepsilon$, respectively.
\end{proof}
\begin{Lem}\label{mu4}
Let $0<\varepsilon < 1/4$. Then, there hold that
\begin{align*}
& \Big\|
\int^{t}_{0}e^{-(t-t')\partial^3_{x}}
f(\cdot,t') \, dt'
\Big\|_{L_x^{120/(26+51\varepsilon)}L_T^{30/(2-3\varepsilon)}} 
 \lesssim T^{\delta}
\| D^{1/12+\varepsilon}_x
f
\|_{L_x^{10/(7-3\varepsilon)}
L_T^{30/(23-12\varepsilon)}} 
\end{align*}
and
\begin{align*}
& \Big\|
\int^{t}_{0}e^{-(t-t')\partial^3_{x}}
f(\cdot,t') \, dt'
\Big\|_{L_x^{60/(13-12\varepsilon)}L_T^{15/(1-4\varepsilon)}} \nonumber \\
& \qquad \lesssim T^{\delta}
\big(\| f\|_{L_x^{10/(7-8\varepsilon)}L_T^{30/(23-32\varepsilon)}}
 +
 \| D_x^{1/12+\varepsilon}
f\|_{L_x^{10/(7-8\varepsilon)}L_T^{30/(23-32\varepsilon)}} \big)
\end{align*}
\end{Lem}
for some  $\delta>0$.
\begin{proof}
We define now
$$
T^{(1)}_z f (x,t)
:=D_x^{-\frac{2z}{3}-\frac{1}{12}}
\int^t_0 
e^{i(t-t')\partial^3_{x}} f(\cdot,t') \, dt'
$$
and
$$
T^{(2)}_z f (x,t)
:=(1+D^2_x)^{-\frac{1}{8}z-\frac{1}{24}}
\int^t_0 
e^{i(t-t')\partial^3_{x}} f(\cdot,t') \, dt'.
$$
We can proceed similarly to the proof of Lemma \ref{PrL24L8/3} using
(see \cite[eq. (3.55) and (3.51)]{KPV1993})
\begin{equation*} 
\| T^{(1)}_{i\gamma}f\|_{L^{60/13}_xL^{15}_T}
=
\Big\|
D^{-1/12}_x
\int^t_0 e^{i(t-t')\partial^3_{x}} 
D_x^{-2i\gamma/3}  f(\cdot,t')dt' \Big\|_{L^{60/13}_xL^{15}_T}
\lesssim 
T^{1/6} \|f \|_{L^{10/7}_xL^{30/23}_T},
\end{equation*}
\begin{equation*}
\| T^{(1)}_{1+i\gamma}f\|_{L^{2}_x L^{\infty}_T}
=
\Big\|
D^{-3/4}_x \int^t_0 e^{i(t-t')\partial^3_{x}} D_x^{-2i\gamma/3}f(\cdot,t')dt' 
\Big\|_{L^{2}_x L^{\infty}_T}
\lesssim T^{1/2}
\|f \|_{L^2_xL^2_T} 
\end{equation*}
and the estimate
(see \cite[eq. (3.55) ]{KPV1993} and \cite[Theorem A.]{Sjo2011})
\begin{align*} 
\| T^{(2)}_{i\gamma}f\|_{L^{60/13}_xL^{15}_T}
& =
\Big\|
(1+D^2_x)^{-\frac{1}{24}}\int^t_0 e^{i(t-t')\partial^3_{x}} 
(1+D^2_x)^{-\frac{i\gamma}{8}}  f(\cdot,t')dt' \Big\|_{L^{60/13}_xL^{15}_T}
\\
& \lesssim 
T^{1/6} \|f \|_{L^{10/7}_xL^{30/23}_T},
\end{align*}
and
\begin{align*}
\| T^{(2)}_{1+i\gamma}f\|_{L^{6}_x L^{\infty}_T}
& =
\Big\|
(1+D^2_x)^{-\frac{1}{6}} \int^t_0 e^{i(t-t')\partial^3_{x}} (1+D^2_x)^{-\frac{i\gamma}{8}}f(\cdot,t')dt' 
\Big\|_{L^{6}_x L^{\infty}_T}\\
& \lesssim T^{1/2}
\|f \|_{L^2_xL^2_T}. 
\qedhere
\end{align*}
\end{proof}

\subsection{Fonseca-Linares-Ponce pointwise formula}

In this section, we present a pointwise formula derived in \cite{FLP2015} that enables the commutation of fractional powers $|x|^r$ and the Airy semigroup $e^{-t\partial^3_{x}}$ with the necessary adjustments. Specifically, for $u_0 \in H^{s}(\R)\cap L^2(|x|^{2r}dx)$, where $0 < s < 2$ and $0 < r \leq s/2$, the following identity holds:
\begin{equation} \label{FLP2}
|x|^{r}e^{-t\partial^3_{x}}u_0(x)=e^{-t\partial^3_{x}}(|x|^{r}u_0)(x)+e^{-t\partial^3_{x}}\{\Psi_{t,r}(\widehat{u}_0)(\xi)\}^{\vee}(x)
\end{equation}
for all $t > 0$ and almost any $x \in \R$. Moreover, the $L^2$--norm of the last term above can be controlled as follows:
\begin{equation} \label{FLP3}
\|\{\Psi_{t,r}(\widehat{u}_0)(\xi)\}^{\vee}\|_{L_x^2}
\lesssim 
(1+t) 
\big(
\|u_0\|_{L^2}
+\|D^{2r}_xu_0\|_{L^2}
\big).
\end{equation}
It is important to note that in \cite{FLP2015}, they specifically considered the case of $s = 2\alpha$ and $r = \alpha$ for $0 < \alpha < 1$. However, a careful analysis of their proof reveals that this result can be extended to the slightly more general situation described here.

\section{Main estimates}
\label{Sect:Keyestimates}

In this section, we gather several interpolation formulas that will prove useful in the subsequent discussions. 
\subsection{Interpolation formulas}
We commence by revisiting a well-known estimate, the proof of which we include here for the sake of completeness.

\begin{Lem}\label{Interpolat_formula}
Let 
$p, p_1, p_2, p_3, q, q_1, q_2, q_3 \in [1,\infty]$ and $\theta_1, \theta_2 \in [0,1]$ such that
\begin{equation*}
\frac{1}{p}
=
\frac{\theta_1}{p_1}
+
\frac{\theta_2}{p_2}
+
\frac{1-\theta_1-\theta_2}{p_3}
\quad \text{and} \quad
\frac{1}{q}
=
\frac{\theta_1}{q_1}
+
\frac{\theta_2}{q_2}
+
\frac{1-\theta_1-\theta_2}{q_3},
\end{equation*}
that is $(1/p,1/q)$ is a point in the convex hull generated by $(1/p_1,1/q_1)$, $(1/p_2,1/q_2)$ and $(1/p_3,1/q_3)$. Then,
\begin{equation*}
\|f\|_{L^p_x L^q_T}
\leq
\|f\|_{L^{p_1}_x L^{q_1}_T}^{\theta_1} \, 
\|f\|_{L^{p_2}_x L^{q_2}_T}^{\theta_2} \,
\|f\|_{L^{p_3}_x L^{q_3}_T}^{1-\theta_1-\theta_2}.
\end{equation*}
\end{Lem}

\begin{proof}
Since
\begin{equation*}
1
=
\frac{1}{\frac{q_1}{q \theta_1}}
+
\frac{1}{\frac{q_2}{q \theta_2}}
+
\frac{1}{\frac{q_3}{q (1-\theta_1-\theta_2)}},
\end{equation*}
the H\"{o}lder inequality allows for the following representation:

\begin{align*}
\|f\|_{L^q_T}^q
& = 
\int_0^T 
|f|^{q \theta_1}
|f|^{q \theta_2}
|f|^{q (1-\theta_1 - \theta_2)}\,dt \\
& = 
\int_0^T 
(|f|^{q_1})^{\frac{q \theta_1}{q_1}}
(|f|^{q_2})^{\frac{q \theta_2}{q_2}}
(|f|^{q_3})^{\frac{q (1-\theta_1 - \theta_2)}{q_3}}
\,dt \\
& \leq 
\big\|
(|f|^{q_1})^{\frac{q \theta_1}{q_1}}
\big\|_{L_T^{\frac{q_1}{q \theta_1}}}
\, 
\big\|
(|f|^{q_2})^{\frac{q \theta_2}{q_2}}
\big\|_{L_T^{\frac{q_2}{q \theta_2}}}
\, 
\big\|
(|f|^{q_3})^{\frac{q (1-\theta_1 - \theta_2)}{q_3}}
\big\|_{L_T^{\frac{q_3}{q (1-\theta_1 - \theta_2)}}} \\
& =
\|f\|_{L_T^{q_1}}^{q \theta_1}
\,
\|f\|_{L_T^{q_2}}^{q \theta_2}
\,
\|f\|_{L_T^{q_2}}^{q (1-\theta_1-\theta_2)}.
\end{align*}
Hence,
\begin{equation}\label{eq:step1}
\|f\|_{L^q_T}
\leq
\|f\|_{L^{q_1}_T}^{\theta_1} \, 
\|f\|_{L^{q_2}_T}^{\theta_2} \,
\|f\|_{L^{q_3}_T}^{1-\theta_1-\theta_2}.
\end{equation}
Similarly, utilizing the relation
\begin{equation*}
1
=
\frac{1}{\frac{p_1}{p \theta_1}}
+
\frac{1}{\frac{p_2}{p \theta_2}}
+
\frac{1}{\frac{p_3}{p (1-\theta_1-\theta_2)}},
\end{equation*}
we get from \eqref{eq:step1} that 
\begin{align*}
\|f\|_{L^p_x L^q_T}^p
&= 
\int_{\R} \|f\|_{L^q_T}^p \, dx 
\leq 
\int_{\R} 
\|f\|_{L^{q_1}_T}^{p \theta_1} \, 
\|f\|_{L^{q_2}_T}^{p \theta_2} \,
\|f\|_{L^{q_3}_T}^{p(1-\theta_1-\theta_2)}
\, dx \\
& =
\int_{\R} 
(\|f\|_{L^{q_1}_T}^{p_1})^{\frac{p \theta_1}{p_1}} \, 
(\|f\|_{L^{q_2}_T}^{p_2})^{\frac{p \theta_2}{p_2}} \,
(\|f\|_{L^{q_3}_T}^{p_3})^{\frac{p(1-\theta_1-\theta_2)}{p_3}}
\, dx \\
& \leq
\Big\| 
(\|f\|_{L^{q_1}_T}^{p_1})^{\frac{p \theta_1}{p_1}}
\Big\|_{L_x^{\frac{p_1}{p \theta_1}}} \,
\Big\| 
(\|f\|_{L^{q_2}_T}^{p_2})^{\frac{p \theta_2}{p_2}}
\Big\|_{L_x^{\frac{p_2}{p \theta_2}}} \,
\Big\| 
(\|f\|_{L^{q_3}_T}^{p_3})^{\frac{p(1-\theta_1-\theta_2)}{p_3}}
\Big\|_{L_x^{\frac{p_3}{p (1-\theta_1-\theta_2)}}} \\
& =
\|f\|_{L^{p_1}_x L^{q_1}_T}^{p\theta_1} \, 
\|f\|_{L^{p_2}_x L^{q_2}_T}^{p\theta_2} \,
\|f\|_{L^{p_3}_x L^{q_3}_T}^{p(1-\theta_1-\theta_2)}. 
\qedhere
\end{align*} 
\end{proof}

Next we state a convenient version for our later analysis of the
Gagliardo--Nirenberg inequality for fractional derivatives.

\begin{Lem}\label{Lem:GaglNiren}
Let $\varepsilon >0$. Then, for 
$\theta:=(1+12\varepsilon)/(13+12\varepsilon)$ we have that
\begin{equation}\label{eqMas1}
\|\partial_xu\|_{L_x^{p}L_T^{q}}
\lesssim \|u\|_{L_x^{p_1}L_T^{q_1}}^{\theta}
\|D^{1/12+\varepsilon}_x \partial_xu\|_{L_x^{\infty}L_T^{2}}^{1-\theta}   ,
\end{equation}
while for 
$\theta:=12/(13+12\varepsilon)$
\begin{equation}\label{eqMas2}
\|D^{1/12+\varepsilon}_xu\|_{L_x^{p}L_T^{q}}
 \lesssim \|u\|_{L_x^{p_1}L_T^{q_1}}^{\theta}\|D^{1/12+\varepsilon}_x \partial_xu\|_{L_x^{\infty}L_T^{2}}^{1-\theta}   ,
\end{equation}
provided that
\begin{equation*}
\frac{1}{p}
= \frac{\theta}{p_1}, \quad
\frac{1}{q}
= \frac{\theta}{q_1} + \frac{1-\theta}{2}.
\end{equation*}
\end{Lem}

\begin{proof}
Firstly, recall from \cite[Lemma 3.3]{MasSeg}
the interpolation formula for fractional derivatives
\begin{equation}\label{MasSegInequality}
 \|D_x^{s}f\|_{L_x^{p}L_T^{q}}
 \lesssim \|D_x^{s_1}f\|_{L_x^{p_1}L_T^{q_1}}^{\theta}
 \|D_x^{s_2}f\|_{L_x^{p_2}L_T^{q_2}}^{1-\theta}   
\end{equation}
for  $p_i, q_i \in (1,\infty)$, $s_i\in \R$ and $\theta\in (0,1)$ related as 
\begin{equation*}
\frac{1}{p}=\frac{\theta}{p_1}+\frac{1-\theta}{p_2},\quad \frac{1}{q}=\frac{\theta}{q_1}+\frac{1-\theta}{q_2},\quad s=\theta s_1+(1-\theta) s_2.    
\end{equation*}
Furthermore, it is allowed to take $p_2=\infty$ in \eqref{MasSegInequality} 
if one uses the interpolation result of O. Blasco for vector valued functions of bounded mean oscillations (\cite[Corollary 1]{Blasco}). Finally, to derive \eqref{eqMas1} and \eqref{eqMas2} we also need to take into account that
$D_x=\mathcal H \partial_x$, $\partial_x=\mathcal H D_x$
and that the Hilbert transform satisfies the following boundedness properties (see \cite[Theorem 2]{Burk} and
\cite[Theorem 1.3]{RubRuiTorr1986})
\begin{equation*}
\| \mathcal H f \|_{L_x^{p}L_T^{q}}
\lesssim \| f \|_{L_x^{p}L_T^{q}}, \qquad
p,q \in (1,\infty)
\end{equation*}
and
\begin{equation*}
\| \mathcal H f \|_{\textsc{BMO}_x L_T^{q}}
\lesssim \| f \|_{L_x^\infty L_T^{q}}, \qquad
q \in (1,\infty). 
\qedhere
\end{equation*}
\end{proof}

 The following lemma is a consequence of the fractional Caffarelli-Kohn-Nirenberg inequality established in \cite{NguSqua2018}.

\begin{Lem}\label{NguSqas}
Let $r, \tau, \gamma > 0$, 
$\sigma \leq 0$,
$0 < a < 1$,
$0 \leq \theta \leq 1$
and $0 \leq \alpha \leq 1/2$
satisfy  
\begin{equation*}
\gamma = a\sigma + (1-a)r \quad \text{and} \quad
\frac{1}{\tau}
+\gamma =a \big( \frac{1-\theta}{2}-\theta \alpha\big)+(1-a) 
\big(\frac{1}{2}+r\big).
\end{equation*}
Then, there exists $\delta > 0$ such that
\begin{equation*}
\| |x|^\gamma f \|_{  L^{6/(\theta\alpha+\theta)}_{T} L_x^\tau}
\lesssim T^\delta \Big(
\|D_x^{\theta\alpha} f\|_{ L^{6/(\theta\alpha+\theta)}_{T} L_x^{2/(1-\theta)}}^a +
\|f\|_{ L^{6/(\theta\alpha+\theta)}_{T} L_x^{2/(1-\theta)}}^a 
\Big) \big\|{|x|^r f} \big\|_{L_T^\infty L^2_x}^{1-a}.
\end{equation*}
\end{Lem}
\begin{proof}
Recall from \cite[Theorem 1.1]{NguSqua2018} that
\begin{equation*}
\| |x|^\gamma f \|_{L_x^\tau}
\lesssim |f|^a_{W_x^{\theta\alpha,2/(1-\theta)}} \,
\big\|{|x|^r f} \big\|_{L^2_x}^{1-a},
\end{equation*}
where
\begin{equation*}
|f|_{W_x^{s,p}}
:=
\big( \int_{\R} \int_{\R} \frac{|f(x)-f(y)|^p}{|x-y|^{1+sp}}\, dx \, dy \big)^{1/p}.
\end{equation*}
Since (see for example \cite[Theorem 5, p.155]{Stein1970})
\begin{equation*}
|f|_{W_x^{s,p}} + \|f\|_{L^p_x} \sim \|D_x^sf\|_{L^p_x}
+ \|f\|_{L^p_x}, \quad p \geq 2,
\end{equation*}
we can write
\begin{equation*}
\| |x|^\gamma f \|_{L_x^\tau} \lesssim
\Big( \|D_x^{\theta\alpha} f\|_{L_x^{2/(1-\theta)}}^a
+ \|f\|_{L_x^{2/(1-\theta)}}^a
\Big) \big\|{|x|^r f} \big\|_{L^2_x}^{1-a}.
\end{equation*}
Next, the H\"older  inequality leads us to
\begin{align*}
 \| |x|^\gamma f \|_{ L_T^{6/(\theta\alpha+\theta)} L_x^\tau}
 &\lesssim \Big( \Big\|
\|D_x^{\theta\alpha} f\|_{L_x^{2/(1-\theta)}}^a
\Big\|_{ L_T^{6/(\theta\alpha+\theta)}} +
\Big\|\|f\|_{L_x^{2/(1-\theta)}}^a \Big\|_{ L_T^{6/(\theta\alpha+\theta)}}
\Big)\big\|{|x|^r f} \big\|_{L^\infty_TL^2_x}^{1-a} \\
&    = 
\Big(\|D_x^{\theta\alpha} f\|_{L_T^{6a/(\theta\alpha+\theta)} L_x^{2/(1-\theta)}}^a
+\|f\|_{L_T^{6a/(\theta\alpha+\theta)} L_x^{2/(1-\theta)}}^a
\Big)
\big\|{|x|^r f} \big\|_{L^\infty_TL^2_x}^{1-a} \\
&   \lesssim 
T^\delta 
\Big( \|D_x^{\theta\alpha} f\|_{L_T^{6/(\theta\alpha+\theta)} L_x^{2/(1-\theta)}}^a
+\|f\|_{L_T^{6/(\theta\alpha+\theta)} L_x^{2/(1-\theta)}}^a \Big)
\big\|{|x|^r f} \big\|_{L^\infty_TL^2_x}^{1-a}
\end{align*}
for some $\delta >0$.
\end{proof}

\subsection{Commutators}
\label{Sect:Commutators}

The following is our first commutator result.
\begin{Lem}\label{Prop:Commutator2}
Let 
$0<r,\rho<1$ such that $r+\rho<1$,
$1 \leq \Tilde{p}_1 < 1/(1-r-\rho)$
and 
$1 \leq q,p \leq \infty$
satisfying
$1/q+1=1/\Tilde{p}_1 + 1/p$.
Then,
\begin{equation*}
\|\{[\langle \xi \rangle^{-\rho},D_\xi^{r}]f\}^{\wedge}  \|_{L^{q}_x}
\lesssim 
\|\widehat{f}\|_{L^{p}_x}. 
\end{equation*}
\end{Lem}

\begin{proof}
Define
\begin{equation*}
K_{\rho}(x)
:= 
\int_{\R}  e^{-i x \xi} \langle \xi \rangle^{-\rho} \, d\xi.
\end{equation*}
It is known that (see for example \cite[Lemma 0.3.9]{Sogge})
\begin{equation}\label{eq:decay}
|K_{\rho}(x)|
\lesssim\frac{1}{|x|^{1-\rho}(1+|x|)^2}.
\end{equation}
By using the estimate 
\begin{align*}
&\big|\{[\langle \xi \rangle^{-\rho},D_\xi^{r}]f\}^{\wedge}(x)\big|
 = 
\big|\{
\langle \xi \rangle^{-\rho}D_\xi^{r}f
-
D_\xi^{r}\big(\langle \xi \rangle^{-\rho}f\big)
\}^{\wedge}(x)\big| \\
& \quad   \sim
\big|
K_{\rho} \ast (|\cdot|^r \widehat{f}) (x)
-
|x|^r 
K_{\rho} \ast \widehat{f} (x)
\big| 
=
\Big|
\int_{\R}
\big(
|y|^r 
-
|x|^r 
\big)
K_{\rho}(x-y) \widehat{f}(y)
\, dy
\Big| \\
& \quad  \leq
\int_{\R}
|x-y|^r \,
|K_{\rho}(x-y)| \, 
|\widehat{f}(y)|
\, dy 
 =
(|\cdot|^r |K_{\rho}|) \ast |\widehat{f}| (x),
\end{align*}
it   follows from the Young  inequality that
\begin{align*}
\|\{[\langle \xi \rangle^{-\rho},D_\xi^{r}]f\}^{\wedge}  \|_{L^{q}_x}
& \lesssim 
\big\||\cdot|^r K_{\rho} \big\|_{L^{\Tilde{p}_1}_x}
\|\widehat{f}\|_{L^{p}_x}.
\end{align*}
Now, inequality  \eqref{eq:decay}
implies
\begin{equation}\label{eq:K2rLp}
\big\||\cdot|^r K_{\rho} \big\|_{L^{\Tilde{p}_1}_x}^{\Tilde{p}_1}
 \lesssim 
\int_0^1 \frac{dx}{|x|^{(1-r-\rho)\Tilde{p}_1}} 
+
\int_1^\infty \frac{dx}{(1+|x|)^{2\Tilde{p}_1}} 
< \infty,
\end{equation}
provided    
\begin{equation*}
1 \leq \Tilde{p}_1 < \frac1{1-r-\rho}.
\qedhere    
\end{equation*}
\end{proof}

Next, we establish a consequence of the Leibniz rule for fractional derivatives.

\begin{Lem} \label{PropComm1} 
If $0<\alpha<1$, 
$1 \leq p \leq 2 \leq q < \infty$ 
and $1/2 = 1/p - 1/q$ we have that
\begin{equation*}
\|[g,D^{\alpha}_x]f\|_{L^2_x}
\lesssim 
\|D_x^{\alpha}g\|_{L^{q}_{x}} 
\, 
\|\widehat{f}\|_{L_x^p}.
\end{equation*}
\end{Lem}

\begin{proof}
For $1/2 = 1/q + 1/p'$, it is known from \cite{KPV1993} 
that
\begin{equation*}
\|D_x^{\alpha}(fg) -gD_x^{\alpha}f -fD_x^{\alpha}g \|_{L^2_x}
 \lesssim \|D_x^{\alpha}g\|_{L^{q}_{x}} \| f\|_{L_x^{p'}}.
\end{equation*} 
So, the Hölder and Hausdorff-Young inequalities imply that
for $1/p+1/p'=1$ and $1 \leq p \leq 2$,
\begin{align*}
\|[g,D^{\alpha}_x]f\|_{L^2_x}
& \leq 
\|D_x^{\alpha}(fg) -gD_x^{\alpha}f -fD_x^{\alpha}g \|_{L^2_x}
+
\|fD_x^{\alpha}g \|_{L^2_x}
\\
& \lesssim  \|D_x^{\alpha}g\|_{L^{q}_{x}} \| f\|_{L_x^{p'}}
\lesssim 
 \|D_x^{\alpha}g\|_{L^{q}_{x}}
 \| \widehat{f}\|_{L_x^{p}}. \qedhere
\end{align*} 
\end{proof}

In Section \ref{Subsubsec:NL1} we will apply the previous result to the function 
$g(\xi)
:=
e^{it\xi^3}
\langle \xi \rangle^{-\rho}$, 
therefore Lemma \ref{lem:Dg} below will be very helpful.





\begin{Lem}\label{lem:Dg}
Let $0<r,\rho<1$ such that 
$r+\rho<1$ and $2r<\rho$.
Then, for any $t>0$ and
$q \in  (\max\{2,1/(\rho-2r)\},\infty)$ one has 
\begin{equation*}
\Big\|
D^{r}_{\xi}\Big(\frac{e^{it\xi^3}}{\langle\xi\rangle^{\rho}}\Big)
\Big\|_{L_{\xi}^{q}}
\lesssim  1+t^\delta
\end{equation*}
for certain $\delta>0$.
\end{Lem}

\begin{proof}
Using properties \eqref{eq:5} and \eqref{eq:6}, we have
\begin{align*}
 \big\|
D^{r}_{\xi}
\big(
e^{it\xi^3}
\langle\xi\rangle^{-\rho}
\big)
\big\|_{L_{\xi}^{q}}
 &\leq 
\big\|
D^{r}_{\xi}
\big(
e^{it\xi^3}
\langle\xi\rangle^{-\rho}
\big)
\big\|_{L_{\xi}^{q}}
+
\big\|
\langle\xi\rangle^{-\rho}
\big\|_{L_{\xi}^{q}} \nonumber \\
& \sim 
\big\|
\mathscr{D}^{r}_{\xi}
\big(
e^{it\xi^3}
\langle\xi\rangle^{-\rho}
\big)
\big\|_{L_{\xi}^{q}}
+
\big\|
\langle\xi\rangle^{-\rho}
\big\|_{L_{\xi}^{q}}
\nonumber \\&
 \leq 
\big\|
\langle\xi\rangle^{-\rho}
\mathscr{D}^{r}_{\xi}
e^{it\xi^3}
\big\|_{L_{\xi}^{q}} 
+
\big\|
\mathscr{D}^{r}_{\xi}
\langle\xi\rangle^{-\rho}
\big\|_{L_{\xi}^{q}} 
+
\big\|
\langle\xi\rangle^{-\rho}
\big\|_{L_{\xi}^{q}}
\nonumber \\
& \sim 
\big\|
\langle\xi\rangle^{-\rho}
\mathscr{D}^{r}_{\xi}
e^{it\xi^3}
\big\|_{L_{\xi}^{q}} 
+
\big\|
D^{r}_{\xi}
\langle\xi\rangle^{-\rho}
\big\|_{L_{\xi}^{q}} 
+
\big\|
\langle\xi\rangle^{-\rho}
\big\|_{L_{\xi}^{q}}.
\end{align*}
The third term on the right-hand side of the above inequality is finite if $q>1/\rho$. As for the second term, it is sufficient to choose
$q > \max\{2,1/(r+\rho)\}$
and utilize the
Hausdorff-Young inequality
in combination with \eqref{eq:K2rLp}.
Finally, by applying \cite[Lemma 2.2]{BusJi2018}, the first term is estimated as follows for some $\delta>0$:
\begin{equation*}
\big\|
\langle\xi\rangle^{-\rho}
\mathscr{D}^{r}_{\xi}
e^{it\xi^3}
\big\|_{L_{\xi}^{q}}
 \lesssim 
t^{\delta}
\big\|
\langle\xi\rangle^{-(\rho-2r)}
\big\|_{L_{\xi}^{q}}
\lesssim 
t^{\delta},
\end{equation*}
provided  $q>1/(\rho-2r)$.
\end{proof} 

\section{Proof of Theorem \ref{Th:Main}}
\label{Sect:ProofTh}

As is standard, our aim is to apply the Banach fixed-point theorem to the mapping
\begin{equation*} 
\Phi(u)
:=e^{-t\partial^3_{x}}u_0
-\int^t_{0}e^{-(t-t')\partial^3_{x}} \partial_x (u^{4})(\cdot,t') \, dt'.
\end{equation*}
Therefore, it suffices to demonstrate that this mapping is a contraction on a suitable subspace of
$L^\infty_T H^{1/12+\varepsilon}_x
\cap L^\infty_T L^2_x(|x|^{2r}dx)$, 
as detailed below. \\

For some $\rho$ and $T>0$, which we will fix later, let us consider the complete metric space
$$
X^{\rho}_T:=\{u \, : \, \|u\|_{X_T} \leq \rho \},
$$
equipped with the norm
\vspace{-\baselineskip}
\begin{equation*}
\|u\|_{X_T}
:= 
\|u\|_{L_T^{\infty}H_x^{1/12+\varepsilon}}
+
\||x|^{r}u\|_{L^{\infty}_{T}L^2_x}
+
\sum_{j=2}^{7} \nu^T_{j}(u)
+
\sum_{j=1}^{6} \mu^T_{j}(u),
\end{equation*}
where
\vspace{-\baselineskip}
\vspace{-\baselineskip}
\begin{multicols}{2}
\begin{align*}    
\nu^T_{2}(u)
& := (1+T)^{-\rho}\|u\|_{L_x^{42/13}L_T^{21/4}}
,\\ 
\nu^T_{3}(u)
& := \|u\|_{L_x^{60/13}L_T^{15}}
,\\ 
\nu^T_{4}(u)
& := T^{-1/6}\|u\|_{L_x^{10/3}L_{T}^{30/7}}, \\
& \\
\nu^T_{5}(u)
& := T^{-1/6}
\|D_x^{1/12} u
\|_{L_x^{10/3}L_{T}^{30/7}},\\
\nu^T_{6}(u)
& := 
\| \partial_x u\|_{L_x^{\infty}L_T^{2}},  \\
\nu^T_{7}(u)
& := 
\|
D_x^{1/12} \partial_x u
\|_{L_x^{\infty}L_{T}^{2}},
\end{align*}  
\end{multicols}
and 
\begin{align*}  
& \mu^T_{1}(u)
:= T^{-1/2}\|u\|_{L_x^{2}L_T^{2}},
&  \mu^T_{4}(u)
:=
\|D_x^{\theta \alpha} u\|_{L_T^{6/(\theta\alpha+\theta)} L_x^{2/(1-\theta)}},  \\
& \mu^T_{2}(u)
 := \|u\|_{L_x^{120/(26+51\varepsilon)}L_T^{30/(2-3\varepsilon)}},
& \mu^T_{5}(u)
 := 
\|
D_x^{1/12+\varepsilon} \partial_x u
\|_{L_x^{\infty}L_{T}^{2}}, \quad \, \,  \\
&\mu^T_{3}(u)
 := \|u\|_{L_T^{6 /(\theta\alpha+\theta)} L_x^{2/(1-\theta)}},
& \mu^T_{6}(u)
 := \|u\|_{L_x^{60/(13-12\varepsilon)}L_T^{15/(1-4\varepsilon)}}. \,
\end{align*}  
Throughout this proof, we fix  
\begin{equation}\label{eq:epsthetaalpha}
\varepsilon:=\frac{633}{5000}, \quad 
R:=\frac{51}{1000}, \quad 
\theta:=\frac{23}{25}, \quad 
\alpha:=\frac{41}{100}.
\end{equation}

\begin{Rem}
To provide a better view of the proof, numerology, and technical constraints, it is essential to note that controlling the norms of the nonlinear contribution of $\Phi(u)$ involves several considerations. Throughout our analysis, we make frequent use of Lemma \ref{Interpolat_formula} to interpolate various $L^p_x L^q_T$-norms of  $u$. This interpolation involves terms such as $\nu_3^T(u)$, $\mu_1^T(u)$, $\mu_2^T(u)$, or $\mu_6^T(u)$.
As a consequence, we cannot arbitrarily reduce the parameter $\varepsilon$ to a very small value. On the other hand, in the treatment of $N\!L_1$ (see Section \ref{Subsubsec:NL1} below), we rely on Lemma \ref{lem:Dg}, which imposes the condition that $q>1/(1/12 + \varepsilon - 2r)$. However, for the application of Lemma \ref{Interpolat_formula}, we must also consider smaller values of $q$. Therefore, we need to restrict the range of $r$ to be within the interval $0 < r < R \ll 1/24 + \varepsilon/2$.
It is also crucial to mention that the parameters $\theta$ and $\alpha$ play a vital role in estimating $N\!L_{3,2}$. See Section \ref{Subsubsec:NL3} for more details.
\end{Rem}
\begin{Rem}
  It is worth noting that
 the $\nu_j^T$ terms were  previously  considered in \cite[p. 585]{KPV1993} for  establishing  the well-posedness in $H^{1/12}(\R)$ in the unweighted case ($r=0$).  
 Therefore, we already know that
\begin{equation}\label{eq:wdnu}
\sum_{j=2}^{7} \nu^T_{j}(\Phi(u))
\lesssim
\|u_0\|_{H^{1/12}}
+
T^{\delta} 
\Big(\max_{j=2, \cdots, 7} \nu^T_{j}(u)\Big)^4
\lesssim 
\|u_0\|_{H^{1/12 + \varepsilon}}
+
T^{\delta} \|u\|_{X_T}^4
\end{equation}
for certain $\delta>0$.
\end{Rem}


\subsection{Analysis of the $L_T^{\infty}H_x^{1/12+\varepsilon}$--norm of $\Phi(u)$} 
\label{Subsec:nu1}
The Plancherel formula readily yields
$$
\| e^{-t\partial^3_{x}}u_0\|_{L_T^{\infty}H_x^{1/12+\varepsilon}}
\lesssim \| u_0\|_{H_x^{1/12+\varepsilon}}.
$$
To control the nonlinear term, observe that  \eqref{nah3.3} and the chain rule \eqref{Sta4.6} guarantee 
\begin{align}\label{eq:25above}
& \Big\|
\int^{t}_{0}e^{-(t-t')\partial^3_{x}}\partial_x( u^4) \, dt'
\Big\|_{L_T^{\infty}H_x^{1/12+\varepsilon}} \nonumber \\
& \qquad 
\lesssim
\Big\| \partial_x
\int^{t}_{0}e^{-(t-t')\partial^3_{x}}
D^{1/12+\varepsilon}_x( u^4) \, dt'
\Big\|_{L_T^{\infty}L_x^{2}} 
+ 
\Big\| 
\partial_x
\int^{t}_{0}e^{-(t-t')\partial^3_{x}}
 u^4 \, dt'
\Big\|_{L_T^{\infty}L_x^{2}} \nonumber
\\
& \qquad \lesssim \|D^{1/12+\varepsilon}_x( u^4)  \|_{L_x^{1}L_T^{2}}+\| u^4  \|_{L_x^{1}L_T^{2}}
\nonumber  \\
& \qquad \lesssim  \|u\|_{L_x^{3r_1}L_T^{3s_1}}^3
\|D^{1/12+\varepsilon}_xu\|_{L_x^{p}L_T^{q}} 
+
T^{1/10} \| u \|_{L_x^{4}L_T^{10}}^4,
\end{align}
provided that
\begin{equation}
\label{eq:chain1}
1
=
\frac{1}{r_1} + \frac{1}{p}, \quad
\frac{1}{2}
=
\frac{1}{s_1} + \frac{1}{q}.
\end{equation}
On the other hand, estimate \eqref{eqMas2} implies
\begin{equation*}
\|D^{1/12+\varepsilon}_xu\|_{L_x^{p}L_T^{q}}
\lesssim 
\|u\|_{L_x^{p_1}L_T^{q_1}}^{\theta_0}
\|D^{1/12+\varepsilon}_x \partial_x u\|_{L_x^{\infty}L_T^{2}}^{1-\theta_0}
\end{equation*}
for 
\begin{equation}
\label{eq:theta0p1q1}
\theta_0:=\frac{12}{13+12\varepsilon}, \quad
\frac{1}{p}
= \frac{\theta_0}{p_1}, \quad
\frac{1}{q}
= \frac{\theta_0}{q_1} + \frac{1-\theta_0}{2}.
\end{equation}
Now, if we set 
$r_1:=127/100$,
$s_1:=287/100$
and take $p,q, p_1$, $q_1$ satisfying \eqref{eq:chain1}
and \eqref{eq:theta0p1q1}, 
it is possible to find by Lemma \ref{Interpolat_formula} the numbers
$\theta_1, \theta_1', \theta_1'', 
\theta_2, \theta_2', \theta_2'' \in (0,1)$ such that
\begin{align*}
\|u\|_{L_x^{3 r_1}L_T^{3 s_1}}
& \lesssim 
T^{\theta_2/2}
(\nu^T_{3}(u))^{\theta_1}\, (\mu^T_{1}(u))^{\theta_2} \, (\mu^T_{2}(u))^{1-\theta_1-\theta_2}, \\
\|u\|_{L_x^{p_1}L_T^{q_1}}
& \lesssim 
T^{\theta_2'/2}
(\nu^T_{3}(u))^{\theta_1'}\, (\mu^T_{1}(u))^{\theta_2'} \, (\mu^T_{2}(u))^{1-\theta_1'-\theta_2'}, \\
\|u\|_{L_x^{4}L_T^{10}}
& \lesssim 
T^{\theta_2''/2}
(\nu^T_{3}(u))^{\theta_1''}\, (\mu^T_{1}(u))^{\theta_2''} \, (\mu^T_{2}(u))^{1-\theta_1''-\theta_2''}.
\end{align*}

Consequently, all the above inequalities lead us to the following bound:
\begin{equation}\label{eq:wdH1/12eps}
\|\Phi(u)\|_{L_T^{\infty}H_x^{1/12+\varepsilon}}
\lesssim 
\|u_0\|_{H^{1/12 + \varepsilon}}
+
T^{\delta} \|u\|_{X_T}^4 
\end{equation}
for some $\delta>0$.

\subsection{Analysis of $\mu_1^T(\Phi(u))$} 
\label{Subsec:mu2}
Inequalities \eqref{2.5lem2.4} and \eqref{eq:Sobolevinclusion} yield
$$
\| e^{-t\partial^3_{x}}u_0\|_{
L_T^{2}
L_x^{2}}
\lesssim 
T^{1/2} \|u_0\|_{H^{1/12 + \varepsilon}}.
$$
Meanwhile, utilizing   estimate \eqref{eq:GinebreTsutsumi1} and the H\"{o}lder  inequality, we obtain
\begin{align}\label{eq:above29}
& \Big\|
\int^{t}_{0}e^{-(t-t')\partial^3_{x}}\partial_x( u^4) \, dt'
\Big\|_{L_x^{2}L_T^{2}}  
 \lesssim 
 T^{1/2}\Big\|
\int^{t}_{0}e^{-(t-t')\partial^3_{x}}\partial_x( u^4) \, dt'
\Big\|_{L_T^{\infty}L_x^{2}}  \nonumber \\
 & \qquad \lesssim
T^{1/2}\|u^3\partial_xu\|_{L_T^{8/7}L_x^{8/7}} 
 \lesssim
 T^{1/2}\|u^3\|_{L_x^{8/7}L_T^{8/3}} \|\partial_xu\|_{L_x^{\infty}L_T^{2}}
\nonumber \\
 & \qquad=
T^{1/2}
\|u\|_{L_x^{24/7}L_T^{8}}^3 \nu_6^T(u).
\end{align}
Moreover, an application Lemma \ref{Interpolat_formula} reveals that
\begin{equation*}
\|u\|_{L_x^{24/7}L_T^{8}}
\lesssim 
T^{\theta_2/2}
(\nu^T_{3}(u))^{\theta_1}\, (\mu^T_{1}(u))^{\theta_2} \, (\mu^T_{2}(u))^{1-\theta_1-\theta_2}
\end{equation*}
for suitable 
$\theta_1, \theta_2 \in (0,1)$.
Hence, we can find $\delta>0$ that satisfies
\begin{equation}\label{eq:wdmu1}
\mu_1^T(\Phi(u))
\lesssim 
\|u_0\|_{H^{1/12 + \varepsilon}}
+
T^{\delta} \|u\|_{X_T}^4.
\end{equation}

\subsection{Analysis of $\mu_{2}^T(\Phi(u))$} 
\label{Subsec:m4}
Lemma \ref{PrL24L8/3} immediately leads to
$$
\| e^{-t\partial^3_{x}}u_0\|_{L_x^{120/(26+51\varepsilon)}L_T^{30/(2-3\varepsilon)}}
\lesssim 
 \|u_0\|_{H^{1/12+\varepsilon}}.
$$
Applying Lemma \ref{mu4} and utilizing the Leibniz rule for fractional derivatives, along with the inequality \eqref{Sta4.6}, we conclude that
\begin{align*}
& \Big\|
\int^{t}_{0}e^{-(t-t')\partial^3_{x}}\partial_x( u^4) \, dt'
\Big\|_{L_x^{120/(26+51\varepsilon)}L_T^{30/(2-3\varepsilon)}} \\
& \qquad \lesssim T^{\delta}
\| D^{1/12+\varepsilon}_x(
u^3 \partial_xu)\|_{L_x^{10/(7-3\varepsilon)}L_T^{30/(23-12\varepsilon)}}
 \\
& \qquad \lesssim T^{\delta}
\big(
\| u^3 D^{1/12+\varepsilon}_x\partial_xu\|_{L_x^{10/(7-3\varepsilon)}L_T^{30/(23-12\varepsilon)}}
+\| \partial_xu\|_{L_x^{p_2}L_T^{q_2}}\| D^{1/12+\varepsilon}_x(u^3) \|_{L_x^{p_1}L_T^{q_1}} 
\big)
\\
& \qquad \lesssim 
T^{\delta}
\big(
\| u\|_{L_x^{30/(7-3\varepsilon)}L_T^{45/(4-6\varepsilon)}}^3\|D^{1/12+\varepsilon}_x\partial_xu\|_{L_x^{\infty}L_T^{2}}\\
& \qquad \qquad \quad
+\| \partial_xu\|_{L_x^{p_2}L_T^{q_2}} \| u^2 \|_{L_x^{r_1}L_T^{s_1}} \| D^{1/12+\varepsilon}_xu\|_{L_x^{p}L_T^{q}}
\big),
\end{align*}
where 
\begin{equation}\label{eq:mu2_p1q1p2q2}
 \frac{7-3\varepsilon}{10}=\frac{1}{p_1}+\frac{1}{p_2},\quad \frac{23-12\varepsilon}{30}=\frac{1}{q_1}+\frac{1}{q_2}    
\end{equation}
and
\begin{equation}\label{eq:mu2_r1s1pq}
 \frac{1}{p_1}=\frac{1}{r_1}+\frac{1}{p},\quad \frac{1}{q_1}=\frac{1}{s_1}+\frac{1}{q}.   
\end{equation}
  On the other hand, we obtain from  \eqref{eqMas1} that
\begin{equation*}
\|\partial_xu\|_{L_x^{p_2}L_T^{q_2}}\lesssim\|u\|_{L_x^{p_3}L_T^{q_3}}^{\theta_0}\|D^{1/12+\varepsilon}_x \partial_xu\|_{L_x^{\infty}L_T^{2}}^{1-\theta_0}   
\end{equation*}
where
\begin{equation}
\label{eq:mu2_theta'0p2q2}
\theta_0:=\frac{1+12\varepsilon}{13+12\varepsilon}, \quad
\frac{1}{p_2}
= \frac{\theta_0}{p_3}, \quad
\frac{1}{q_2}
= \frac{\theta_0}{q_3} + \frac{1-\theta_0}{2}.
\end{equation}
We also have from  \eqref{eqMas2} and H\"older's inequality
\begin{equation*}
\|D^{1/12+\varepsilon}_xu\|_{L_x^{p}L_T^{q}}
\lesssim
\|u\|_{L_x^{p_4}L_T^{q_4}}^{\theta_0'}
\|D^{1/12+\varepsilon}_x \partial_xu\|_{L_x^{\infty}L_T^{2}}^{1-\theta_0'} 
\leq 
T^{\delta}
\|u\|_{L_x^{p_4}L_T^{9q_4/8}}^{\theta_0'}
(\mu^T_{5}(u))^{1-\theta_0'} ,
\end{equation*}
where 
\begin{equation}
\label{eq:mu2_theta0p4q4}
\theta_0':=\frac{12}{13+12\varepsilon}, \quad
\frac{1}{p}
= \frac{\theta_0'}{p_4}, \quad
\frac{1}{q}
= \frac{\theta_0'}{q_4} + \frac{1-\theta_0'}{2}.
\end{equation}

Taking
$r_1:=23/10$,
$s_1:=15/2$,
$p_2:=24$,
$q_2:=233/100$ 
and  $p,q$, $p_1,q_1$, $p_3,q_3$, $p_4,q_4$,  satisfying  
\eqref{eq:mu2_p1q1p2q2},
\eqref{eq:mu2_r1s1pq},
\eqref{eq:mu2_theta'0p2q2}
and
\eqref{eq:mu2_theta0p4q4},  by Lemma \ref{Interpolat_formula}
we find  
$\theta_1, \theta_1', \theta_1'', \theta_1'''$, $\theta_2, \theta_2', \theta_2'' , \theta_2''' \in (0,1)$ such that 
\begin{align*}
\|u\|_{L_x^{30/(7-3\varepsilon)}L_T^{45/(4-6\varepsilon)}}&  \lesssim T^{\theta_2/2}
(\nu^T_{3}(u))^{\theta_1}\, (\mu^T_{1}(u))^{\theta_2} \, (\mu^T_{2}(u))^{1-\theta_1-\theta_2},
\\
\|u\|_{L_x^{2 r_1}L_T^{2 s_1}}
 & \lesssim  T^{\theta_2'/2}
(\nu^T_{3}(u))^{\theta_1'}\, (\mu^T_{1}(u))^{\theta_2'} \, (\mu^T_{2}(u))^{1-\theta_1'-\theta_2'},
\\
\|u\|_{L_x^{p_3}L_T^{q_3}}
& \lesssim T^{\theta_2''/2}
(\nu^T_{3}(u))^{\theta_1''}\, (\mu^T_{1}(u))^{\theta_2''} \, (\mu^T_{2}(u))^{1-\theta_1''-\theta_2''},
\\
\|u\|_{L_x^{p_4}L_T^{9q_4/8}}
& \lesssim T^{\theta_2'''/2}
(\nu^T_{3}(u))^{\theta_1'''}\, (\mu^T_{1}(u))^{\theta_2'''} \, (\mu^T_{2}(u))^{1-\theta_1'''-\theta_2'''}.
\end{align*}

Combining all the above, we conclude that
\begin{equation}\label{eq:wdmu2}
\mu_{2}^T(\Phi(u))
\lesssim 
\|u_0\|_{H^{1/12 + \varepsilon}}
+
T^{\delta} \|u\|_{X_T}^4 
\end{equation}
for some $\delta>0$.

\subsection{Analysis of $\mu_3^T(\Phi(u))$} 
\label{Subsec:mu5}

By H\"older's inequality and \eqref{eq:KPV1989} (taken with $\alpha=0$) we obtain
\begin{align*}
& 
\|e^{-t\partial^3_{x}}u_0\|_{L_T^{6 /(\theta\alpha+\theta)} L_x^{2/(1-\theta)}}
\leq 
T^\delta
\|e^{-t\partial^3_{x}}u_0\|_{L_T^{6 /\theta} L_x^{2/(1-\theta)}}
\lesssim
T^\delta
 \|u_0\|_{L^2}
\lesssim 
T^\delta
 \|u_0\|_{H^{1/12 + \varepsilon}}.
\end{align*} 
Moreover, inequality  \eqref{eq:GinebreTsutsumi1} with $p_2=q_2=8$  gives
\begin{align*}
\Big\| 
\int^{t}_{0}e^{-(t-t')\partial^3_{x}} \partial_x(u^4) \, dt' \Big\|_{L_T^{6 /(\theta\alpha+\theta)} L_x^{2/(1-\theta)}} 
& \lesssim
T^{\delta} \Big\| 
\int^{t}_{0}e^{-(t-t')\partial^3_{x}} \partial_x(u^4) \, dt' \Big\|_{L_T^{6 /\theta}L_x^{2/(1-\theta)}}
\\
&  \lesssim
T^{\delta}\|u^3\partial_xu\|_{L_T^{8/7}L_x^{8/7}},
\end{align*}
and using the same argument as in \eqref{eq:above29}
 we can conclude 
\begin{equation}\label{eq:wdmu3}
\mu_3^T(\Phi(u))
\lesssim 
\|u_0\|_{H^{1/12 + \varepsilon}}
+
T^{\delta} \|u\|_{X_T}^4 
\end{equation}
for some $\delta>0$.

\subsection{Analysis of $\mu_4^T(\Phi(u))$} 
\label{Subsec:mu6}
By \eqref{eq:KPV1989} and taking $\theta\alpha/2<1/12+\varepsilon$ we have 
\begin{align*}
\|D_x^{\theta\alpha}e^{-t\partial^3_{x}}u_0\|_{L_T^{6/(\theta\alpha+\theta)} L_x^{2/(1-\theta)}}
& \lesssim 
\|D_x^{\theta\alpha/2}
e^{-t\partial^3_{x}}
D^{\theta\alpha/2}_xu_0\|_{L_T^{6/(\theta\alpha+\theta)} L_x^{2/(1-\theta)}}
\\
&  \lesssim
\|D^{\theta\alpha/2}_x u_0\|_{L_x^2}
\lesssim
\|u_0\|_{H^{1/12+\varepsilon}}.
\end{align*}
Inequality \eqref{eq:KPV1989}, the Minkowski    inequality  
(since 
$6/(6-\theta\alpha-\theta)
>
2/(1+\theta)$
for the particular choice $\theta:=23/25$ and
$\alpha:=41/100$), and the H\"older  inequality imply
\begin{align*}
& \Big\| D_x^{\theta\alpha}
\int^{t}_{0}e^{-(t-t')\partial^3_{x}} \partial_x(u^4) \, dt' \Big\|_{L_T^{6/(\theta\alpha+\theta)} L_x^{2/(1-\theta)}}
\lesssim 
\|u^3\partial_xu\|_{L_T^{6/(6-\theta\alpha-\theta)} L_x^{2/(1+\theta)} }
\\
& \qquad 
\lesssim 
\|u^3\partial_xu\|_{L_x^{2/(1+\theta)} L_T^{6/(6-\theta\alpha-\theta)} }
\lesssim \|u\|_{ L_x^{3r_1}L_T^{3s_1}}^3 \,\|\partial_xu\|_{ L_x^{p}L_T^{q}}
\end{align*}
for
\begin{equation}
\label{MU4_rspq}
\frac{1+\theta}{2}
=
\frac{1}{r_1} + \frac{1}{p}, \quad
\frac{6-\theta\alpha-\theta}{6}
=
\frac{1}{s_1} + \frac{1}{q}.
\end{equation}
Next, estimate \eqref{eqMas1} produces
\begin{equation*}
\|\partial_xu\|_{ L_x^{p}L_T^{q}} 
\lesssim 
\|u\|_{L_x^{p_1}L_T^{q_1}}^{\theta_0}
\|D^{1/12+\varepsilon}_x \partial_x u\|_{L_x^{\infty}L_T^{2}}^{1-\theta_0}
\end{equation*}
for 
\begin{equation}
\label{MU4_theta'0p1q1}
\theta_0:=\frac{1+12\varepsilon}{13+12\varepsilon}, \quad
\frac{1}{p}
= \frac{\theta_0}{p_1}, \quad
\frac{1}{q}
= \frac{\theta_0}{q_1} + \frac{1-\theta_0}{2}.
\end{equation}
Finally, if we set 
$r_1:=1111/1000$,
$s_1:=299/100$
and take $p,q, p_1$, $q_1$ satisfying \eqref{MU4_rspq}
and \eqref{MU4_theta'0p1q1}, by Lemma \ref{Interpolat_formula}
it is possible to find 
$\theta_1, \theta_1', 
\theta_2, \theta_2'\in (0,1)$ such that
\begin{align*}
\|u\|_{ L_x^{3r_1}L_T^{3s_1}}
& \lesssim 
T^{\theta_2/2}
(\nu^T_{3}(u))^{\theta_1}\, (\mu^T_{1}(u))^{\theta_2} \, (\mu^T_{2}(u))^{1-\theta_1-\theta_2}, \\
\|u\|_{L_x^{p_1}L_T^{q_1}}
& \lesssim 
T^{\theta_2'/2}
(\nu^T_{3}(u))^{\theta_1'}\, (\mu^T_{1}(u))^{\theta_2'} \, (\mu^T_{2}(u))^{1-\theta_1'-\theta_2'}.
\end{align*}

Hence, there exists $\delta>0$ such that
\begin{equation}\label{eq:wdmu4}
\mu_4^T(\Phi(u))\lesssim 
\|u_0\|_{H^{1/12 + \varepsilon}}
+
T^{\delta} \|u\|_{X_T}^4.
\end{equation}

\subsection{Analysis of $\mu_{5}^T(\Phi(u))$} 
\label{Subsec:mu9}
Inequality  \eqref{ineq36} leads  to
$$
\|
\partial_x
e^{-t\partial^3_{x}}
D^{1/12+\varepsilon}_x
u_0\|_{L_x^{\infty}L_T^{2}}
\lesssim 
\| D^{1/12+\varepsilon}_xu_0\|_{L_x^2}
\lesssim \|u_0\|_{H^{1/12+\varepsilon}}.
$$

Moreover, by \eqref{kpv-es-3} 
\begin{align*}
& \Big\| \partial^2_x
\int^{t}_{0}e^{-(t-t')\partial^3_{x}} D^{1/12+\varepsilon}_x(u^4) \, dt'
\Big\|_{L_x^{\infty}L_T^{2}} 
\lesssim
\|D^{1/12+\varepsilon}_x(u^4)\|_{L_x^{1}L_T^{2}},
\end{align*}
and proceeding as in \eqref{eq:25above} above we can get 
\begin{equation}\label{eq:wdmu5}
\mu_5^T(\Phi(u))
\lesssim 
\|u_0\|_{H^{1/12 + \varepsilon}}
+
T^{\delta} \|u\|_{X_T}^4
\end{equation}
for certain $\delta>0$.

\subsection{Analysis of $\mu_6^T(\Phi(u))$} 
\label{Subsec:mu6}

Lemma \ref{PrL24L8/3} inmediatly leads to
$$
\| e^{-t\partial^3_{x}}u_0\|_{L_x^{60/(13-12\varepsilon)}L_T^{15/(1-4\varepsilon)}}
\lesssim 
 \|u_0\|_{H^{1/12+\varepsilon}}.
$$
Applying Lemma \ref{mu4}, we have 
\begin{align*}
& \Big\|
\int^{t}_{0}e^{-(t-t')\partial^3_{x}}\partial_x( u^4) \, dt'
\Big\|_{L_x^{60/(13-12\varepsilon)}L_T^{15/(1-4\varepsilon)}} \\
& \qquad \lesssim T^{\delta} 
\big( \| u^3 \partial_xu\|_{L_x^{10/(7-8\varepsilon)}L_T^{30/(23-32\varepsilon)}}+\| D^{1/12+\varepsilon}_x(
u^3 \partial_xu)\|_{L_x^{10/(7-8\varepsilon)}L_T^{30/(23-32\varepsilon)}}\big).
\end{align*}
Next, the H\"older   inequality yields to
\begin{align*}
 &\|u^3 \partial_x u \|_{L_x^{10/(7-8\varepsilon)}L_T^{30/(23-32\varepsilon)}}
\\
& \qquad \lesssim  
\|u^3\|_{L_x^{10/(7-8\varepsilon)}L_T^{15/(4-16\varepsilon)}}\|\partial_xu\|_{L_x^{\infty}L_T^{2}}
\lesssim 
\|u\|_{L_x^{30/(7-8\varepsilon)}L_T^{45/(4-16\varepsilon)}}^3 \,\nu^T_{6}(u). 
\end{align*}

Furthermore, utilizing the Leibniz rule for fractional derivatives, along with the inequality \eqref{Sta4.6}, we conclude that
\begin{align*}
& \| D^{1/12+\varepsilon}_x(
u^3 \partial_xu)\|_{L_x^{10/(7-8\varepsilon)}L_T^{30/(23-32\varepsilon)}}
 \\
& \, \lesssim
\| u^3 D^{1/12+\varepsilon}_x\partial_xu\|_{L_x^{10/(7-8\varepsilon)}L_T^{30/(23-32\varepsilon)}}
+\| \partial_xu\|_{L_x^{p_2}L_T^{q_2}}\| D^{1/12+\varepsilon}_x(u^3) \|_{L_x^{p_1}L_T^{q_1}} 
\\
& \, \lesssim 
\| u\|_{L_x^{30/(7-8\varepsilon)}L_T^{45/(4-16\varepsilon)}}^3\|D^{1/12+\varepsilon}_x\partial_xu\|_{L_x^{\infty}L_T^{2}}
+\| \partial_xu\|_{L_x^{p_2}L_T^{q_2}} \| u \|_{L_x^{2r_1}L_T^{2s_1}}^2 \| D^{1/12+\varepsilon}_xu\|_{L_x^{p}L_T^{q}},
\end{align*}
provided that 
\begin{equation}\label{eq:mu6_p1q1p2q2}
 \frac{7-8\varepsilon}{10}=\frac{1}{p_1}+\frac{1}{p_2},\quad \frac{23-32\varepsilon}{30}=\frac{1}{q_1}+\frac{1}{q_2}    
\end{equation}
and
\begin{equation}\label{eq:mu6_r1s1pq}
 \frac{1}{p_1}=\frac{1}{r_1}+\frac{1}{p},\quad \frac{1}{q_1}=\frac{1}{s_1}+\frac{1}{q}.   
\end{equation}
  On the other hand, we obtain from  \eqref{eqMas1}
  and the H\"older  inequality that
\begin{equation*}
\|\partial_xu\|_{L_x^{p_2}L_T^{q_2}}\lesssim\|u\|_{L_x^{p_3}L_T^{q_3}}^{\theta_0}\|D^{1/12+\varepsilon}_x \partial_xu\|_{L_x^{\infty}L_T^{2}}^{1-\theta_0} \lesssim T^{\delta}\|u\|_{L_x^{p_3}L_T^{2q_3}}^{\theta_0}\|D^{1/12+\varepsilon}_x \partial_xu\|_{L_x^{\infty}L_T^{2}}^{1-\theta_0},   
\end{equation*}
where
\begin{equation}
\label{eq:mu6_theta'0p2q2}
\theta_0:=\frac{1+12\varepsilon}{13+12\varepsilon}, \quad
\frac{1}{p_2}
= \frac{\theta_0}{p_3}, \quad
\frac{1}{q_2}
= \frac{\theta_0}{q_3} + \frac{1-\theta_0}{2}.
\end{equation}
We also have from  \eqref{eqMas2}  that
\begin{equation*}
\|D^{1/12+\varepsilon}_xu\|_{L_x^{p}L_T^{q}}
\lesssim
\|u\|_{L_x^{p_4}L_T^{q_4}}^{\theta_0'}
\|D^{1/12+\varepsilon}_x \partial_xu\|_{L_x^{\infty}L_T^{2}}^{1-\theta_0'} 
\end{equation*}
for 
\begin{equation}
\label{eq:mu6_theta0p4q4}
\theta_0':=\frac{12}{13+12\varepsilon}, \quad
\frac{1}{p}
= \frac{\theta_0'}{p_4}, \quad
\frac{1}{q}
= \frac{\theta_0'}{q_4} + \frac{1-\theta_0'}{2}.
\end{equation}

Taking
$r_1:=63/25$,
$s_1:=27/2$,
$p_2:=27$,
$q_2:=23/10$ 
and  $p,q$, $p_1,q_1$, $p_3,q_3$, $p_4,q_4$,  satisfying  
\eqref{eq:mu6_p1q1p2q2},
\eqref{eq:mu6_r1s1pq},
\eqref{eq:mu6_theta'0p2q2}
and
\eqref{eq:mu6_theta0p4q4},  by Lemma \ref{Interpolat_formula}
we find  
$\theta_1, \theta_1', \theta_1'', \theta_1'''$, $\theta_2, \theta_2', \theta_2'' , \theta_2''' \in (0,1)$ such that 
\begin{align*}
\|u\|_{L_x^{30/(7-8\varepsilon)}L_T^{45/(4-16\varepsilon)}}&  
\lesssim 
(\nu^T_{3}(u))^{\theta_1}\, (\mu^T_{2}(u))^{\theta_2} \, (\mu^T_{6}(u))^{1-\theta_1-\theta_2},
\\
\|u\|_{L_x^{2 r_1}L_T^{2 s_1}}
 & \lesssim  
(\nu^T_{3}(u))^{\theta_1'}\, (\mu^T_{2}(u))^{\theta_2'} \, (\mu^T_{6}(u))^{1-\theta_1'-\theta_2'},
\\
\|u\|_{L_x^{p_3}L_T^{2q_3}}
& \lesssim  
(\nu^T_{3}(u))^{\theta_1''}\, (\mu^T_{2}(u))^{\theta_2''} \, (\mu^T_{6}(u))^{1-\theta_1''-\theta_2''},
\\
\|u\|_{L_x^{p_4}L_T^{q_4}}
& \lesssim  
(\nu^T_{3}(u))^{\theta_1'''}\, (\mu^T_{2}(u))^{\theta_2'''} \, (\mu^T_{6}(u))^{1-\theta_1'''-\theta_2'''}.
\end{align*}

Combining all the above, we conclude that
\begin{equation}\label{eq:wdmu6}
\mu_{6}^T(\Phi(u))
\lesssim 
\|u_0\|_{H^{1/12 + \varepsilon}}
+
T^{\delta} \|u\|_{X_T}^4.
\end{equation}

\subsection{Analysis of the $L^{\infty}_{T}L^2_x(|x|^{2r})$--norm of $\Phi(u)$}
\label{Subsec:weightednorm}

First, we note that
\begin{align*}
\||x|^{r}\Phi(u)\|_{L^{\infty}_{T}L^2_x}
&\leq 
\||x|^{r}e^{-t\partial^3_{x}}u_0\|_{L^{\infty}_{T} L^2_x}
+
\Big\|
|x|^{r}
\int^t_{0}e^{-(t-t')\partial^3_{x}} \partial_x (u^4)\,dt'
\Big\|_{L^{\infty}_{T}L^2_x} 
\nonumber\\ &
=: 
L + N\!L. 
\end{align*}

The linear term $L$ can be controlled by \eqref{FLP2}, \eqref{FLP3} and thePlancherel  formula \eqref{2.5lem2.4}
\begin{align}\label{esu0}
L
&\leq  
\|e^{-t\partial^3_{x}}(|x|^{r}u_0)\|_{L^{\infty}_{T}L^2_x}
+
\|e^{-t\partial^3_{x}}(\{\Psi_{t,r}(\widehat{u}_0)(\xi)\}^{\vee})\|_{L^{\infty}_{T}L^2_x}\nonumber\\ 
& \lesssim
\||x|^{r}u_0\|_{L^2}
+
(1+T)
\big(\|u_0\|_{L^2}
+\|D^{2r}_x u_0\|_{L^2}\big)
\lesssim
\||x|^{r}u_0\|_{L^2}
+
\|u_0\|_{H^{1/12 + \varepsilon}},
\end{align}
where we  used \eqref{eq:Sobolevinclusion} and the fact $T \lesssim 1$.\\

To treat  the nonlinear factor $N\!L$, we decompose it as follows:
\begin{align*}
N\!L
& \sim
\Big\|
D_\xi^{r}\Big\{\int^t_{0}e^{i(t-t')\xi^3}
\big(
\partial_x (u^4)
\big)^{\wedge}(\xi,t') \,dt'
\Big\}
\Big\|_{L^{\infty}_{T}L^2_\xi} \\
& =
\Big\|
\int^t_{0}
D_\xi^{r}\Big\{
\frac{e^{i(t-t')\xi^3}}
{\langle \xi \rangle^{1/12+\varepsilon}}
\langle \xi \rangle^{1/12+\varepsilon}
\big(
\partial_x (u^4)
\big)^{\wedge}(\xi,t')
\Big\}
\, dt'
\Big\|_{L^{\infty}_{T}L^2_\xi} \\
& \leq
\Big\|
\int^t_{0}
\Big[
\frac{e^{i(t-t')\xi^3}}
{\langle \xi \rangle^{1/12+\varepsilon}}
,
D_\xi^{r}
\Big]
\Big\{
\langle \xi \rangle^{1/12+\varepsilon}
\big(
\partial_x (u^4)
\big)^{\wedge}(\xi,t')
\Big\}
\, dt'
\Big\|_{L^{\infty}_{T}L^2_\xi} \\
& \quad +
\Big\|
\int^t_{0}
e^{i(t-t')\xi^3}
[\langle \xi \rangle^{-1/12-\varepsilon}
,
D_\xi^{r}]
\Big\{
\langle \xi \rangle^{1/12+\varepsilon}
\big(
\partial_x (u^4)
\big)^{\wedge}(\xi,t')
\Big\}
\, dt'
\Big\|_{L^{\infty}_{T}L^2_\xi} \\
& \quad +
\Big\|
\int^t_{0}
e^{i(t-t')\xi^3}
D_\xi^{r}\Big\{
\big(
\partial_x (u^4)
\big)^{\wedge}(\xi,t')
\Big\}
\, dt'
\Big\|_{L^{\infty}_{T}L^2_\xi}
=: 
N\!L_1
+
N\!L_2
+
N\!L_3.
\end{align*}

We investigate each of the above terms in the next subsections.

\subsubsection{Commutator estimate for $N\!L_1$} 
\label{Subsubsec:NL1}
Let $r \in (0,R)$, for 
$R:=51/1000$. 
By Lemma \ref{PropComm1}, 
Lemma \ref{lem:Dg},
and the Minkowski inequality we deduce 
\begin{align*}
&N\!L_1  \leq 
\int^T_{0}
 \Big\|
\Big[ \frac{e^{i(t-t')\xi^3}} {\langle \xi \rangle^{1/12+\varepsilon}},
D_\xi^{r} \Big]\Big\{\langle \xi \rangle^{1/12+\varepsilon}
\big(\partial_x (u^4)
\big)^{\wedge}(\xi,t')
\Big\}
\Big\|_{L^{\infty}_{T}L^2_\xi}
\, dt'
\\
& \qquad\lesssim    
\int^T_{0}
\Big\|
D_{\xi}^r 
\Big(
\frac{e^{i(t-t')\xi^3}} {\langle \xi \rangle^{1/12+\varepsilon}}
\Big)
\Big\|_{L^{\infty}_{T}L^{q}_\xi }
\Big\|
\Big(\langle \xi \rangle^{1/12+\varepsilon}
\big(\partial_x (u^4)
\big)^{\wedge}(\cdot,t')\Big)^{\wedge}(x)
\Big\|_{L^{p}_x} 
\, dt'
\\ & \qquad 
\lesssim
(1+T^{\delta})
\big(
 \|\partial_x (u^4) \|_{L^{1}_{T}L^{p}_x} 
+ 
 \|D^{1/12+\varepsilon}_x \partial_x (u^4) \|_{L^{1}_{T}L^{p}_x}  
 \big)\\
& \qquad  \lesssim 
T^{\delta}
\big(
 \|u^3 \partial_x u \|_{L^{p}_x L^{p}_{T}}
+ 
 \|D^{1/12+\varepsilon}_x (u^3 \partial_x u ) \|_{L^{p}_x L^{p}_{T}}  
 \big)
\end{align*}
for some $\delta>0$ and $1 \leq p \leq 2$ such that
\begin{equation}\label{eq:pqProp3.5Lem3.6}
\frac{1}{2}
=
\frac{1}{p} - \frac{1}{q} \qquad \text{and} \qquad q>\frac{1}{\frac{1}{12} + \varepsilon -2r}.
\end{equation}
The H\"older   inequality yields to
\begin{equation*}
 \|u^3 \partial_x u \|_{L^{p}_x L^{p}_{T}} \lesssim  
\|u^3\|_{L_x^{p_1}L_T^{q_1}}\|\partial_xu\|_{L_x^{\infty}L_T^{2}}
\lesssim 
\|u\|_{L_x^{3p_1}L_T^{3q_1}}^3 \,\nu^T_6(u)
\end{equation*}
provided that 
\begin{equation}\label{eq:relationsppp1q12}
\frac{1}{p}
= \frac{1}{p_1}, \quad
\frac{1}{p}
= \frac{1}{q_1} + \frac{1}{2}.
\end{equation}
Furthermore,   inequalities \eqref{thA8} and \eqref{Sta4.6} give
\begin{align*}
& \|D^{1/12+\varepsilon}_x (u^3 \partial_x u ) \|_{L^{p}_x L^{p}_{T}} \\
& \qquad
\lesssim 
\|u^3\|_{L_x^{p_1}L_T^{q_1}}\|D^{1/12+\varepsilon}_x\partial_xu\|_{L_x^{\infty}L_T^{2}} 
+\|\partial_xu\|_{L_x^{p_2}L_T^{q_2}}  \|D^{1/12+\varepsilon}_x(u^3)\|_{L_x^{p_3}L_T^{q_3}} 
\\
& \qquad
\lesssim 
\|u\|_{L_x^{3p_1}L_T^{3q_1}}^3 \, \mu_{5}^T(u) +\|\partial_xu\|_{L_x^{p_2}L_T^{q_2}}  \|u^2\|_{L_x^{p_4}L_T^{q_4}}\|D^{1/12+\varepsilon}_xu\|_{L_x^{p_5}L_T^{q_5}}
\\
& \qquad
\lesssim 
\|u\|_{L_x^{3p_1}L_T^{3q_1}}^3 \, \mu_{5}^T(u)
+
\|\partial_xu\|_{L_x^{p_2}L_T^{q_2}} \|u\|_{L_x^{2p_4}L_T^{2q_4}}^2 \|D^{1/12+\varepsilon}_xu\|_{L_x^{p_5}L_T^{q_5}},
\end{align*}
provided that 
\begin{equation}\label{eq:relationsppp5q5}
\frac{1}{p}
= \frac{1}{p_2} + \frac{1}{p_3}, \quad
\frac{1}{p}
= \frac{1}{q_2} + \frac{1}{q_3}, \quad
\frac{1}{p_3}
= \frac{1}{p_4} + \frac{1}{p_5}, \quad
\frac{1}{q_3}
= \frac{1}{q_4} + \frac{1}{q_5}.
\end{equation}
Next, we employ Lemma 
\ref{Lem:GaglNiren}  to get
\begin{equation*}
\|\partial_xu\|_{L_x^{p_2}L_T^{q_2}} 
\lesssim \|u\|_{L_x^{r_1}L_T^{s_1}}^{\theta_0}
\|D^{1/12+\varepsilon}_x \partial_xu\|_{L_x^{\infty}L_T^{2}}^{1-\theta_0}\lesssim \|u\|_{L_x^{r_1}L_T^{s_1}}^{\theta_0}
\, (\mu^T_5(u))^{1-\theta_0}  
\end{equation*}
for 
\begin{equation}
\label{NL1_theta0r1s1}
\theta_0:=\frac{1+12\varepsilon}{13+12\varepsilon}, \quad
\frac{1}{p_2}
= \frac{\theta_0}{r_1}, \quad
\frac{1}{q_2}
= \frac{\theta_0}{s_1} + \frac{1-\theta_0}{2} 
\end{equation}
and
\begin{equation*}
\|D^{1/12+\varepsilon}_xu\|_{L_x^{p_5}L_T^{q_5}}
\lesssim \|u\|_{L_x^{r_2}L_T^{s_2}}^{\theta'_0}\|D^{1/12+\varepsilon}_x \partial_xu\|_{L_x^{\infty}L_T^{2}}^{1-\theta'_0}
\lesssim T^\delta
\|u\|_{L_x^{r_2}L_T^{s_2}}^{\theta'_0}\, (\mu^T_5(u))^{1-\theta'_0},
\end{equation*}
provided that
\begin{equation}
\label{eq:NL1_theta0p5q5}
\theta'_0:=\frac{12}{13+12\varepsilon}, \quad
\frac{1}{p_5}
= \frac{\theta'_0}{r_2}, \quad
\frac{1}{q_5}
= \frac{\theta'_0}{s_2} + \frac{1-\theta'_0}{2}.
\end{equation}

Take $q:=15000/1619$, 
$p_4:= 3 p_1/2$, 
$q_4:=3q_1/2$,  
$r_2:=49/10$,
$s_2:=55/2$  
and 
$p$, $p_1$, $q_1$,   $p_3,q_3$,     $p_5,q_5$, $r_1$, $s_1$,   satisfying  \eqref{eq:pqProp3.5Lem3.6}-\eqref{eq:NL1_theta0p5q5}. Then,  we can apply Lemma \ref{Interpolat_formula} to obtain
$\theta_1, \theta_1', \theta_1''$,  $\theta_2, \theta_2', \theta_2''  \in (0,1)$ such that 
\begin{align*}
\|u\|_{L_x^{3p_1}L_T^{3q_1}}
 & \lesssim  
(\nu^T_{3}(u))^{\theta_1}\, (\mu^T_{2}(u))^{\theta_2} \, (\mu^T_{6}(u))^{1-\theta_1-\theta_2}, 
\\
\|u\|_{L_x^{r_1}L_T^{s_1}}
&  \lesssim  
(\nu^T_{3}(u))^{\theta_1'}\, (\mu^T_{2}(u))^{\theta_2'} \, (\mu^T_{6}(u))^{1-\theta_1'-\theta_2'},
\\
\|u\|_{L_x^{r_2}L_T^{s_2}}
 & \lesssim  
(\nu^T_{3}(u))^{\theta_1''}\, (\mu^T_{2}(u))^{\theta_2''} \, (\mu^T_{6}(u))^{1-\theta_1''-\theta_2'''}.
\end{align*}

Summarizing, we get
\begin{equation}\label{est:NL1}
N\!L_1
\lesssim 
T^{\delta} \|u\|_{X_T}^4
\end{equation}
for some $\delta>0$.

\subsubsection{Commutator estimate for $N\!L_2$} 
\label{Subsubsec:NL2}
Let $r \in (0, 1/24+\varepsilon/2]$ and for simplicity, denote
\begin{equation*}
f(\xi,t'):= 
\langle \xi \rangle^{1/12+\varepsilon}
\big(\partial_x(u^4)\big)^{\wedge}(\xi,t').  
\end{equation*}
By using property \eqref{eq:GinebreTsutsumi1}, Lemma \ref{Prop:Commutator2}, and the
H\"older    and Minkowski  inequalities we obtain
\begin{align}\label{eq:step1NL2}
 N\!L_2
& \sim 
\Big\|
\int^t_{0}e^{i(t-t')\partial_x^3} 
\Big([\langle \xi \rangle^{-1/12-\varepsilon},D_\xi^{r}]f\Big)^{\bigwedge} (x,t')\, dt'
\Big\|_{L^{\infty}_{T}L^2_x} \nonumber \\
& \lesssim
\|\{[\langle \xi \rangle^{-1/12-\varepsilon},D_\xi^{r}]f\}^{\wedge}  \|_{L^{q'_2}_TL^{p'_2}_x} 
\lesssim
\|\widehat{f}\|_{L^{q'_2}_TL^{p}_x} \nonumber\\ 
& \lesssim  
\| \partial_x(u^4)\|_{L^{q'_2}_TL^{p}_x} 
 +\|D_x^{1/12+\varepsilon} \partial_x(u^4)\|_{L^{q'_2}_TL^{p}_x} \nonumber \\ 
& \lesssim 
T^{\delta}
\big(
\|u^3\partial_xu\|_{L^{p}_xL^{q'_2+d}_T} 
+
\|D^{1/12+\varepsilon}_x(u^3\partial_xu)\|_{L^{p}_xL^{q'_2+d}_T}
\big),
\end{align}
provided that 
\begin{equation}\label{eq:GinebYoung}
p_2 \geq 2, \quad
\frac{1}{q_2}
= \frac{1}{6}
- \frac{1}{3 p_2}, \quad
\frac{1}{p_2} + \frac{1}{p_2'}
=
1
= \frac{1}{q_2} + \frac{1}{q_2'},
\end{equation}
\begin{equation*}
\frac{1}{p_2'} + 1
= \frac{1}{\tilde{p_1}} + \frac{1}{p}, \quad
1 \leq \tilde{p_1} \leq \frac{1}{\frac{11}{12}-r-\varepsilon}
\end{equation*}
and $q'_2+d\geq p$.\\

The first term in \eqref{eq:step1NL2}
can be controlled via the H\"older  inequality and 
\eqref{eqMas1} by
\begin{align*}
\|u^3\partial_xu\|_{L^{p}_xL^{q'_2+d}_T} 
& \leq 
\|u^3\|_{L^{r_0}_xL^{s_0}_T} \|\partial_xu\|_{L^{r_1}_xL^{s_1}_T} \\
& \lesssim 
T^\delta 
\|u\|_{L^{3r_0}_xL^{4s_0}_T}^3
\|u\|_{L_x^{r_2}L_T^{s_2}}^{\theta'_0}
\|D^{1/12+\varepsilon}_x \partial_xu\|_{L_x^{\infty}L_T^{2}}^{1-\theta'_0},
\end{align*}
for
\begin{equation}\label{eq:NL2_rsr1s1}
\frac{1}{p}=\frac{1}{r_0}+\frac{1}{r_1}, \quad 
\frac{1}{q'_2+d}=\frac{1}{s_0}+\frac{1}{s_1}
\end{equation}
and
\begin{equation}
\label{eq:NL2_theta'0r1s1r2s2}
\theta'_0:=\frac{1+12\varepsilon}{13+12\varepsilon}, \quad
\frac{1}{r_1}
= \frac{\theta'_0}{r_2}, \quad
\frac{1}{s_1}
= \frac{\theta'_0}{s_2} + \frac{1-\theta'_0}{2}.
\end{equation}

As for the second term in \eqref{eq:step1NL2} we invoke 
\eqref{thA8} and \eqref{Sta4.6},
\begin{align}\label{eq:D1/12epsetc}
& \|D^{1/12+\varepsilon}_x(u^3\partial_xu)\|_{L^{p}_xL^{q'_2+d}_T}
 \nonumber \\
 & \qquad \lesssim 
\|u^3\|_{L_x^{p}L_T^{q_3}}\|D^{1/12+\varepsilon}_x\partial_xu\|_{L_x^{\infty}L_T^{2}} +\|\partial_xuD^{1/12+\varepsilon}_x(u^3)\|_{L^{p}_xL^{q'_2+d}_T} \nonumber \\
& \qquad \lesssim  \|u\|_{L_x^{3p}L_T^{3q_3}}^3 \, \mu_{5}^T(u) +\|\partial_xu\|_{L_x^{p_4}L_T^{q_4}} \|D^{1/12+\varepsilon}_x(u^3)\|_{L_x^{p_5}L_T^{q_5}}  \nonumber\\
& \qquad \lesssim  
 \|u\|_{L_x^{3p}L_T^{3q_3}}^3 \, \mu_{5}^T(u)+\|\partial_xu\|_{L_x^{p_4}L_T^{q_4}} \|u\|_{L_x^{2r_3}L_T^{2s_3}}^2 \|D^{1/12+\varepsilon}_xu\|_{L_x^{p_7}L_T^{q_7}},
\end{align}
where
\begin{equation}\label{eq:NL2_p4q4p5q5}
 \frac{1}{q'_2+d}=\frac{1}{q_3}+\frac{1}{2}, \quad \frac{1}{p}=\frac{1}{p_4}+\frac{1}{p_5},\quad \frac{1}{q'_2+d}=\frac{1}{q_4}+\frac{1}{q_5}   
\end{equation}
and
\begin{equation}\label{eq:NL2_r3s3p7q7}
 \frac{1}{p_5}=\frac{1}{r_3}+\frac{1}{p_7},\quad \frac{1}{q_5}=\frac{1}{s_3}+\frac{1}{q_7}.   
\end{equation}

Furthermore, we can bound some of the factors in \eqref{eq:D1/12epsetc} by means of \eqref{eqMas1} and \eqref{eqMas2} as follow
\begin{equation*}
\|\partial_xu\|_{L_x^{p_4}L_T^{q_4}}\lesssim\|u\|_{L_x^{p_6}L_T^{q_6}}^{\theta'_0}\|D^{1/12+\varepsilon}_x \partial_xu\|_{L_x^{\infty}L_T^{2}}^{1-\theta'_0}   ,
\end{equation*}
and
\begin{equation*}
\|D^{1/12+\varepsilon}_xu\|_{L_x^{p_7}L_T^{q_7}}\lesssim\|u\|_{L_x^{p_8}L_T^{q_8}}^{\theta_0}\|D^{1/12+\varepsilon}_x \partial_xu\|_{L_x^{\infty}L_T^{2}}^{1-\theta_0} 
\end{equation*}
where 
\begin{equation}
\label{eq:NL2_theta'0p6q6}
\frac{1}{p_4}
= \frac{\theta'_0}{p_6}, \quad
\frac{1}{q_4}
= \frac{\theta'_0}{q_6} + \frac{1-\theta'_0}{2}
\end{equation}
and  
\begin{equation}
\label{eq:NL2_theta0p8q8}
\theta_0:=\frac{12}{13+12\varepsilon}, \quad
\frac{1}{p_7}
= \frac{\theta_0}{p_8}, \quad
\frac{1}{q_7}
= \frac{\theta_0}{q_8} + \frac{1-\theta_0}{2}.
\end{equation}

If we set $p_2:=233/50,\,\tilde{p_1}:=1 $,  then by \eqref{eq:GinebYoung} one can easily compute $p, q'_2$.
And we choose $d:=17/100$ such that $q'_2+d\geq p$.
Also taking  $r_0:=7/5$, $s_0:=16/5$,
$r_3:=19/10$, $s_3:=26/5$, $p_4:=17$, $q_4:=11/5 $
and $r_i, s_i$ (i=1,2), 
 $p_j,q_j (j=5,..,8)$,  satisfying 
 \eqref{eq:NL2_rsr1s1}, \eqref{eq:NL2_theta'0r1s1r2s2}, \eqref{eq:NL2_p4q4p5q5}, 
 \eqref{eq:NL2_r3s3p7q7},
 \eqref{eq:NL2_theta'0p6q6}
 and 
\eqref{eq:NL2_theta0p8q8}, by Lemma \ref{Interpolat_formula}  we find  
$\theta_1, \theta_1', \theta_1''$, $ \theta_2, \theta_2', \theta_2''$, $ \theta_3, \theta_3', \theta_3''$, $\theta_4, \theta_4', \theta_4'' \in (0,1)$ such that 
\begin{align*}
\|u\|_{L^{3r_0}_xL^{4s_0}_T} &  \lesssim T^{\theta_2/2}
(\nu^T_{3}(u))^{\theta_1}\, (\mu^T_{1}(u))^{\theta_2} \, (\mu^T_{2}(u))^{1-\theta_1-\theta_2},
\\
\|u\|_{L_x^{r_2}L_T^{s_2}}
& \lesssim  T^{\theta_2'/2}
(\nu^T_{3}(u))^{\theta_1'}\, (\mu^T_{1}(u))^{\theta_2'} \, (\mu^T_{2}(u))^{1-\theta_1'-\theta_2'},
\\
\|u\|_{L_x^{3p}L_T^{3q_3}}
&  \lesssim T^{\theta_2/2}
(\nu^T_{3}(u))^{\theta''_1}\, (\mu^T_{1}(u))^{\theta''_2} \, (\mu^T_{2}(u))^{1-\theta''_1-\theta''_2},
\\
\|u\|_{L_x^{2 r_3}L_T^{2 s_3}}
 & \lesssim  T^{\theta_2'/2}
(\nu^T_{3}(u))^{\theta_3}\, (\mu^T_{1}(u))^{\theta_4} \, (\mu^T_{2}(u))^{1-\theta_3-\theta_4},
\\
\|u\|_{L_x^{p_6}L_T^{q_6}}
& \lesssim T^{\theta_3'/2}
(\nu^T_{3}(u))^{\theta_3'}\, (\mu^T_{1}(u))^{\theta_4'} \, (\mu^T_{2}(u))^{1-\theta_3'-\theta_3'},
\\
\|u\|_{L_x^{p_8}L_T^{q_8}}
& \lesssim T^{\theta_3''/2}
(\nu^T_{3}(u))^{\theta_3''}\, (\mu^T_{1}(u))^{\theta_4''} \, (\mu^T_{2}(u))^{1-\theta_3''-\theta_4''}.
\end{align*}

Combining all the above, we conclude that
\begin{equation}\label{est:NL2}
N\!L_2
\lesssim 
T^{\delta} \|u\|_{X_T}^4
\end{equation}
for some $\delta>0$.

\subsubsection{Boundedness of $N\!L_3$} 
\label{Subsubsec:NL3}

First we observe that by the Plancherel theorem that
\begin{align*}
N\!L_3
& \sim 
\Big\|
\int^t_{0}e^{-(t-t')\partial^3_{x}}
\big\{
|x|^{r}
\partial_x (u^4) (x,t')
\big\}dt'
\Big\|_{L^{\infty}_{T}L^2_x}. 
\end{align*}

By choosing $\varphi \in C_0^{\infty}(\R)$ such that $\varphi =1$ when $|x|<1/2$, while $\varphi =0$ when $|x|\geq 1$, we have
\begin{align*}  
|x|^{r} \partial_x (u^4)
& =
\varphi(x)|x|^{r} \partial_x (u^4)
+
\partial_x\big((1-\varphi(x)|x|^{r}u^{4}\big) 
-
\partial_x\big((1-\varphi(x))|x|^{r}\big)u^{4}.
\end{align*}
Hence, 
\begin{align*}
N\!L_3 &  
\lesssim 
\Big\|\int^{t}_{0} e^{-(t-t')\partial^3_{x}}
\big\{ \varphi(x)|x|^{r} \partial_x (u^4) \big\}
dt' \Big\|_{L^{\infty}_{T}L^2_x}  
\\  & \qquad  
+ 
\Big\| \int^{t}_{0} e^{-(t-t')\partial^3_{x}}
\partial_x\big((1-\varphi(x)|x|^{r}u^{4}\big) 
dt' \Big\|_{L^{\infty}_{T}L^2_x}  
\\  & \qquad  
+
\Big\| \int^t_{0}e^{-(t-t')\partial^3_{x}}
\big\{ 
\partial_x\big((1-\varphi(x))|x|^{r}\big)u^{4}
\big\}
dt' \Big\|_{L^{\infty}_{T}L^2_x}
\\  & 
=:N\!L_{3,1}+N\!L_{3,2}+N\!L_{3,3}.
\end{align*}

\underline{Analysis of $N\!L_{3,1}$}. 
Inequality \eqref{eq:GinebreTsutsumi1} with $p_2=q_2=8$ leads  to the same estimate as in \eqref{eq:above29} above: 
\begin{align*}
N\!L_{3,1}  
&  
 \lesssim  
 \|\varphi(x)|x|^{r}
\partial_x (u^4)
 \|_{L_T^{8/7}L_x^{8/7}}    
\lesssim
\|u^3 \partial_x u\|_{L_T^{8/7}L_x^{8/7}} \lesssim
\|u^3\partial_xu\|_{L_x^{8/7}L_T^{8/7}}.   
\end{align*}

Thus, repeating the argument for the nonlinear part of $\mu^T_1(\Phi(u))$,
 we conclude 
\begin{equation}\label{NL31}
N\!L_{3,1} \lesssim 
T^{\delta} \|u\|_{X_T}^4 
\end{equation}
for some $\delta>0$.\\

\underline{Analysis of $N\!L_{3,3}$}. By using \eqref{eq:GinebreTsutsumi1} with $q_2=90/11,\,p_2=15/2$ and $q'_2=90/79,\,p'_2=15/13$ and then the Minkowski  inequality, we obtain 
\begin{align*}
N\!L_{3,3}
&   
\lesssim  
\|u^4 \|_{L^{90/79}_{T}L^{15/13}_x } 
\lesssim  
{T}^{\delta} \|u^4 \|_{L^{15/4}_{T}L^{15/13}_x } 
\lesssim 
{T}^{\delta}\|u\|^4_{ L^{15}_{T}L^{60/13}_x}  
\leq
{T}^{\delta}(\nu^T_{3}(u))^4,
\end{align*}
since
\begin{equation*}
\big\|
\partial_x\big(|x|^{r}(1-\varphi(x)\big)
\big\|_{L^\infty_x}
\lesssim 1, \quad 0<r<1.
\end{equation*}
It follows that
\begin{equation}\label{NL33}
N\!L_{3,3}
\lesssim 
T^{\delta} \|u\|_{X_T}^4 
\end{equation}
for some $\delta>0$.\\

\underline{Analysis of $N\!L_{3,2}$}.
An application of \eqref{nah3.3}, the
H\"older inequality  and the Minkowski inequality 
produce
\begin{align*}
N\!L_{3,2}
&  
=
\Big\| \partial_x \int^{t}_{0} e^{-(t-t')\partial^3_{x}}
\big((1-\varphi(x)|x|^{r}u^{4}\big) 
dt' \Big\|_{L^{\infty}_{T}L^2_x} 
\lesssim  
\||x|^{r}  u^4 \|_{L^{1}_xL^2_{T}} 
 \\
& \lesssim 
\||x|^{r}u^{c}\|_{L_x^{p_1} L^{q_1}_{T}} \|u^{b}\|_{L_x^{p_2} L^{q_2}_{T}}
 \lesssim 
\| |x|^{r/c} u \|_{L_x^{cp_1} L^{cq_1}_{T}}^{c} \|u\|^{b}_{L_x^{bp_2} L^{bq_2}_{T}} \\
& \lesssim 
\| |x|^{r/c} u \|_{ L^{cq_1}_{T}L_x^{cp_1}}^{c}
\|u\|^{b}_{L_x^{bp_2} L^{bq_2}_{T}}
\end{align*}
with $p_1\geq q_1$ and
\begin{equation}\label{eq:b+c}
4=b+c, \quad 
1=\frac{1}{p_1}+\frac{1}{p_2}, \quad
\frac12=\frac{1}{q_1}+\frac{1}{q_2}.
\end{equation}
Next, we apply Lemma \ref{NguSqas}
\begin{align*}
\| |x|^{r/c} u \|_{ L^{cq_1}_{T}L_x^{cp_1}}
& \lesssim
T^\delta \Big( \|D_x^{\theta\alpha} u\|_{L_T^{6/(\theta\alpha+\theta)} L_x^{2/(1-\theta)}}^a + \|u\|_{L_T^{6/(\theta\alpha+\theta)} L_x^{2/(1-\theta)}}^a \Big)
\big\|{|x|^ru} \big\|_{L_T^\infty L^2_x}^{1-a}
\\ 
&  = T^{\delta} 
\Big((\mu^T_3(u))^{a}
+(\mu^T_4(u))^{a} \Big)\big\|{|x|^ru} \big\|_{L_T^\infty L^2_x}^{1-a}
\end{align*}
for
\begin{equation}\label{eq:NgSqaCond}
cq_1=\frac{6}{\theta\alpha+\theta}, \quad
\frac{r}{c} = a\sigma + (1-a)r, \quad
\frac{1}{c p_1}
+ \frac{r}{c} =a \big( \frac{1-\theta}{2}-\theta \alpha\big)+(1-a) 
\Big(\frac{1}{2}+r\Big).
\end{equation}

Recall the choices of $\varepsilon$, $\theta$ and $\alpha$ given by \eqref{eq:epsthetaalpha}.
For any $r \in (0, 1/24+\varepsilon/2]$, by taking $a:=69/200$ and $c:=153/100$ one can easily compute  $p_1$, $q_1$, $\sigma$ from \eqref{eq:NgSqaCond} satisfying  $p_1\geq q_1$ and  $\sigma<0$ . 
Furthermore, for $b$, $p_2$, $q_2$ verifying  \eqref{eq:b+c},  Lemma \ref{Interpolat_formula}   gives
\begin{equation*}
\|u\|_{L_x^{bp_2} L^{bq_2}_{T}}
\lesssim 
T^{\theta_2/2}(\nu^T_3(u))^{\theta_1} (\mu^T_1(u))^{\theta_2}(\mu^T_2(u))^{(1-\theta_1-\theta_2)},
\end{equation*}
for certain $\theta_1, \theta_2 \in (0,1)$. Hence,
\begin{equation}\label{NL32}
N\!L_{3,2}
\lesssim 
T^{\delta} \|u\|_{X_T}^4 
\end{equation}
for some $\delta>0$.
 
\subsection{Contraction mapping}
\label{Subsec:contraction}
Let $u \in X_T^\rho$. Collecting
\eqref{eq:wdnu},
\eqref{eq:wdH1/12eps},
\eqref{eq:wdmu1},
\eqref{eq:wdmu2},
\eqref{eq:wdmu3},
\eqref{eq:wdmu4},
\eqref{eq:wdmu5},
\eqref{eq:wdmu6},
\eqref{esu0},
\eqref{est:NL1},
\eqref{est:NL2},
\eqref{NL31},
\eqref{NL33}
and 
\eqref{NL32}
we get
\begin{align*}
\|\Phi(u)\|_{X_T}
& \leq 
C 
\big(
\|u_0\|_{H^{1/12 + \varepsilon}}
+
\||x|^{r}u_0\|_{L^2}
\big)
+
C T^{\delta} \rho^4,
\end{align*}
for some positive constants $C$ and $\delta$.
Thus, by taking
\begin{equation*}
\rho:= 2
\big(
\|u_0\|_{H^{1/12 + \varepsilon}}
+
\||x|^{r}u_0\|_{L^2}
\big)
\end{equation*}
and choosing $T>0$ sufficiently small such that
\begin{equation}\label{conditforT}
\frac{\rho}{2} + C T^{\delta} \rho^4 
\leq 
\rho,
\end{equation}
we can justify that $\Phi(u) \in X_T^\rho$, that is 
$\Phi : X_T^\rho \longmapsto X_T^\rho$ is a well-defined mapping. \\

Furthermore, for $u,v \in X_T^\rho$ we observe that
\begin{equation*} 
\Phi(u)-\Phi(v)
=
\int_0^t 
e^{-(t-t') \partial^3_x} 
\partial_x (u^4-v^4)\, dt',
\end{equation*}
and since
\begin{align*}
\partial_x (u^4-v^4)
& = 
\partial_x 
\big(
(u-v)(u^3+u^2v+uv^2+v^3)
\big) \\
& = 
\partial_x(u-v) \, (u^3+u^2v+uv^2+v^3)   \\
& \quad +
(u-v)
\big(
3u^2 \partial_x u
+
2uv \partial_x u
+
u^2 \partial_x v
+
v^2 \partial_x u
+
2 uv \partial_x v
+
3v^2 \partial_x v \big),
\end{align*}
one can essentially proceed as above to deduce that
\begin{align*}
\|\Phi(u) - \Phi(v)\|_{X_T}
& \leq 
C' T^{\delta} \rho^3
\|u- v\|_{X_T},
\end{align*}
for certain $C'>0$. Therefore, if we pick $T>0$ satisfying simultaneously $C' T^{\delta} \rho^3<1$ and \eqref{conditforT} we can guarantee that $\Phi$ is a contraction.

\addtocontents{toc}{\vspace{1\baselineskip}}

\begin{thebibliography}{10}

\bibitem{BenPan}
{\sc A.~Benedek and R.~Panzone}, {\em The space {$L^{p}$}, with mixed norm},
  Duke Math. J., 28 (1961), pp.~301--324.

\bibitem{BeKa}
{\sc Y.~A. Berezin and V.~I. Karpman}, {\em Nonlinear evolution of disturbances
  in plasma and other dispersive media}, Soviet Physics JETP, 24 (1967),
  pp.~1049--1056.

\bibitem{Blasco}
{\sc O.~Blasco}, {\em Interpolation between ${H}^1_{{B}_0}$ and
  ${L}^p_{{B}_1}$}, Stud. Math., 9 (1989), pp.~205--210.

\bibitem{Bour93II}
{\sc J.~Bourgain}, {\em Fourier transform restriction phenomena for certain
  lattice subsets and applications to nonlinear evolution equations. {II}.
  {T}he {K}d{V}-equation}, Geom. Funct. Anal., 3 (1993), pp.~209--262.

\bibitem{Burk}
{\sc D.~L. Burkholder}, {\em A geometric condition that implies the existence
  of certain singular integrals of {B}anach-space-valued functions}, in
  Conference on harmonic analysis in honor of {A}ntoni {Z}ygmund, {V}ol. {I},
  {II} ({C}hicago, {I}ll., 1981), Wadsworth Math. Ser., 1983, pp.~270--286.

\bibitem{BusJi2018}
{\sc E.~Bustamante, J.~Jiménez, and J.~Mejía}, {\em A note on the {O}strovsky
  equation in weighted {S}obolev spaces}, J. Math. Anal. Appl., 460 (2018),
  pp.~1004--1018.

\bibitem{CN2015}
{\sc X.~Carvajal and W.~Neves}, {\em Persistence property in weighted {S}obolev
  spaces for nonlinear dispersive equations}, Quart. Appl. Math., 73 (2015),
  pp.~493--510.

\bibitem{CJZ2022}
{\sc A.~J. Castro, K.~Jabbarkhanov, and L.~Zhapsarbayeva}, {\em The {N}onlinear
  {S}chr\"{o}dinger-{A}iry equation in weighted {S}obolev spaces}, Nonlinear
  Anal., 223 (2022), p.~113068.

\bibitem{CKSTT}
{\sc J.~Colliander, M.~Keel, G.~Staffilani, H.~Takaoka, and T.~Tao}, {\em Sharp
  global well-posedness for {K}d{V} and modified {K}d{V} on {$\mathbb R$} and
  {$\mathbb T$}}, J. Amer. Math. Soc., 16 (2003), pp.~705--749.

\bibitem{FLP2015}
{\sc G.~Fonseca, F.~Linares, and G.~Ponce}, {\em On persistence properties in
  fractional weighted space}, Proc. Amer. Math. Soc., 143 (2015),
  pp.~5353--5367.

\bibitem{FonPn}
{\sc G.~Fonseca and G.~Ponce}, {\em The {IVP} for the {B}enjamin-{O}no equation
  in weighted {S}obolev spaces}, J. Funct. Anal., 260 (2011), pp.~436--459.

\bibitem{GnTst89}
{\sc J.~Ginibre and Y.~Tsutsumi}, {\em Uniqueness of solutions for the
  generalized {K}orteweg-de {V}ries equation}, SIAM J. Math. Anal., 20 (1989),
  pp.~1388--1425.

\bibitem{Grun2005}
{\sc A.~Gr\"unrock}, {\em A bilinear {A}iry-estimate with application to
  g{K}d{V}-3}, Differ. Integral Equ., 18 (2005), pp.~1333--1339.

\bibitem{Guo}
{\sc Z.~Guo}, {\em Global well-posedness of {K}orteweg-de {V}ries equation in
  ${H}^{-3/4}$}, J. Math. Pur. Appl., 91 (2009), pp.~583--597.

\bibitem{Kato1983}
{\sc T.~Kato}, {\em On the {C}auchy problem for the (generalized) {K}orteweg-de
  {V}ries equation}, in Studies in applied mathematics, vol.~8 of Adv. Math.
  Suppl. Stud., Academic Press, New York, 1983, pp.~93--128.

\bibitem{KPV1989}
{\sc C.~E. Kenig, G.~Ponce, and L.~Vega}, {\em On the (generalized)
  {K}orteweg-de {V}ries equation}, Duke Math. J., 59 (1989), pp.~585--610.

\bibitem{KPVCP}
\leavevmode\vrule height 2pt depth -1.6pt width 23pt, {\em The {C}auchy problem
  for {K}orteveg-de {V}ries in {S}obolev spaces of negative indices}, Duke
  Math. J., 71 (1993), pp.~1--21.

\bibitem{KPV1993}
\leavevmode\vrule height 2pt depth -1.6pt width 23pt, {\em Well-posedness and
  scattering results for the generalized {K}orteweg-de {V}ries equation via the
  contraction principle}, Comm. Pure Appl. Math., 46 (1993), pp.~527--620.

\bibitem{Kis2009}
{\sc N.~Kishimoto}, {\em Well-posedness of the {C}auchy problem for the
  {K}orteweg-de {V}ries equation at the critical regularity}, Diff. and Int.
  Eq., 22 (2009), pp.~447--464.

\bibitem{KdV}
{\sc D.~J. Korteweg and G.~de~Vries}, {\em On the change of form of long waves
  advancing in a rectangular canal, and on a new type of long stationary
  waves}, Philos. Mag. (5), 39 (1895), pp.~422--443.

\bibitem{MarMer}
{\sc Y.~Martel and F.~Merle}, {\em Blow up in finite time and dynamics of blow
  up solutions for the {$L^2$}-critical generalized {K}d{V} equation}, J. Amer.
  Math. Soc., 15 (2002), pp.~617--664.

\bibitem{MasSeg}
{\sc S.~Masaki and J.~I. Segata}, {\em On the well-posedness of the generalized
  {K}orteweg -de {V}ries equation in scale-critical {${\hat{L}}^r$}-space},
  Anal. PDE, 9 (2016), pp.~699--725.

\bibitem{MG2019}
{\sc A.~Mu\~noz{-}Garc\'ia}, {\em Cauchy problem for {K}d{V} equation in
  weighted {S}obolev spaces}, Rev. Fac. Cienc., 8 (2019), pp.~83--102.

\bibitem{Nah2012}
{\sc J.~Nahas}, {\em A decay property of solutions to the k-generalized {K}d{V}
  equation}, Adv. Differ. Equat., 17 (2012), pp.~833--858.

\bibitem{NP2009}
{\sc J.~Nahas and G.~Ponce}, {\em On the persistent properties of solutions to
  semi-linear {S}chr\"{o}dinger equation}, Commun. Partial Differ. Equ., 34
  (2009), pp.~1208--1227.

\bibitem{Nar}
{\sc G.~A. Nariboli}, {\em Nonlinear longitudinal dispersive waves in elastic
  rods}, J. Math. Phys. Sci., 4 (1970), pp.~64--73.

\bibitem{NguSqua2018}
{\sc H.~M. Nguyen and M.~Squassina}, {\em Fractional
  {C}affarelli–{K}ohn–{N}irenberg inequalities}, J. Funct. Anal., 274
  (2018), pp.~2661--2672.

\bibitem{RubRuiTorr1986}
{\sc J.~L. Rubio~de Francia, F.~J. Ruiz, and J.~L. Torrea}, {\em
  {C}alder{\'o}n-{Z}ygmund theory for operator-valued kernels}, Adv. Math., 62
  (1986), pp.~7--48.

\bibitem{Sjo2011}
{\sc P.~Sj{\"o}lin}, {\em Radial functions and maximal operators of
  {S}chr\"{o}dinger type}, Indiana Univ. Math. J., 60 (2011), pp.~143--159.

\bibitem{Sogge}
{\sc C.~D. Sogge}, {\em Fourier {I}ntegrals in {C}lassical {A}nalysis},
  vol.~210, Cambridge University Press, 2017.

\bibitem{Stein1970}
{\sc E.~Stein}, {\em Singular Integrals and Differentiability properties of
  functions}, Princeton University Press, 1970.

\bibitem{Tao2007}
{\sc T.~Tao}, {\em Scattering for the quartic generalised {K}orteweg-de {V}ries
  equation}, J. Differ. Equ., 32 (2007), pp.~623--651.

\bibitem{Li2009}
{\sc L.~Zi-Liang}, {\em Application of higher-order {KdV}-{mKdV} model with
  higher-degree nonlinear terms to gravity waves in atmosphere}, Chin. Phys. B,
  18 (2009), pp.~4074--4082.

\end{thebibliography}

\end{document}